\documentclass[preprint,1p,leqno]{elsarticle}

\usepackage{amssymb}
\usepackage{amsmath}

\journal{Somewhere}
\usepackage{amsmath,amsthm,amsfonts}
\usepackage{amssymb}

\newtheorem{lemma}{Lemma}
\newtheorem{theorem}{Theorem}

\usepackage{xcolor}

\usepackage[T1]{fontenc}
\usepackage{amsmath}

\usepackage[pagewise]{lineno}

\usepackage{bigstrut}
\usepackage{subfig}

\usepackage{geometry}
\usepackage{bm}

\usepackage{graphicx}
\usepackage{amsbsy}
\usepackage{amsfonts}
\usepackage{mathrsfs}

\usepackage{etoolbox}                                       
\apptocmd{\thebibliography}{\setlength{\itemsep}{0pt}}{}{}  
\usepackage[hidelinks]{hyperref}
\begin{document}

\begin{frontmatter}



\title{A decoupled energy-stable mixed finite element method for Poisson-Nernst-Planck-Navier-Stokes equations}


\author{Yanmei Peng}
\ead{pengyanmei_math@outlook.com}

\author{Ying Yang}
\ead{yangying@lsec.cc.ac.cn}

\author[]{Yi Zhang\corref{cor1}}
\ead{zhangyi_aero@hotmail.com}

\cortext[cor1]{Corresponding author}

\affiliation[1]{organization={School of Mathematics and Computational Science, Guilin University of Electronic Technology},
city={Guilin},
country={China}}

\affiliation[2]{organization={Guangxi Colleges and Universities Key Laboratory of Data Analysis and Computation},
city={Guilin},
country={China}}

\affiliation[3]{organization={Guangxi Applied Mathematics Center (GUET)},
city={Guilin},
country={China}}




\begin{abstract}

We propose a novel linearized mixed finite element method for the Poisson-Nernst-Planck-Navier-Stokes (PNPNS) system. Specifically, the method combines a staggered time discretization that eliminates the need for expensive nonlinear solvers by carefully treating nonlinear terms in a time-staggered manner, with a mimetic spatial discretization that preserves the exact structure of the discrete de Rham complex. Both semi-discrete scheme and its fully discrete counterpart are developed, which preserve key physical properties, including conservation of the total mass and energy stability. Under appropriate assumptions on the initial data, a rigorous theoretical analysis is carried out for the fully discrete scheme. Numerical experiments using mimetic spectral elements are presented to demonstrate the properties and to verify the accuracy and effectiveness of the proposed decoupled approach.

\end{abstract}




\end{frontmatter}



\section{Introduction}
\label{sec1}

Electrokinetic phenomena involve complex interplays between fluid dynamics and the transport of charged ions. They are ubiquitous in a wide array of physical, biophysical, and industrial engineering processes. Typical applications range from the modeling of biological ion channels and cellular microenvironments to the design of microfluidic devices, drug delivery systems, and water desalination technologies  \cite{bazant2004diffuse, eisenberg1998ionic, roubivcek2007incompressible,tado2016analysis}. Mathematically, these electro-hydrodynamic processes are most comprehensively described by the coupled Poisson-Nernst-Planck-Navier-Stokes (PNPNS) equations. In this multiphysics system, the Nernst-Planck equations govern the drift-diffusion of different ionic species driven by concentration gradients and electrical fields. The Poisson equation correlates the electrical potential to the distribution of these charge carriers. And the Navier-Stokes equations describe the macroscopic motion of the incompressible fluid, which is in turn subjected to the electrical body forces generated by the ions \cite{lu2010poisson, rubinstein1990electro}.

Due to the highly nonlinear and strongly coupled nature of the physical mechanisms, spanning multiple spatial and temporal scales, the mathematical analysis of the PNPNS system has attracted significant interest. Early foundational works primarily focused on the steady-state properties and the qualitative behavior of the electro-diffusion of ions \cite{ jerome2012analysis, park1997qualitative}. Subsequently, extensive mathematical theories regarding the well-posedness, global existence, uniqueness, and regularity of weak and smooth solutions for the transient system under various boundary conditions have been established \cite{bothe2014global, ryham2006energetic-thesis, schmuck2009analysis, wang2016generalized}. Recent advancements have further extended these analyses to the study of quasi-neutral limits, boundary layer problems, and partial regularities in critical functional spaces \cite{wang2019quasi, gong2021partial}.

While theoretical analyses provide a solid foundation, designing reliable and efficient numerical methods for the PNPNS system remains a formidable challenge. A physically faithful numerical scheme must strictly preserve the intrinsic structural properties of the continuous model at the discrete level. These properties include, for example, the non-negativity of ionic concentrations, the conservation of global mass for each species, and the dissipation of the total free energy. In the broader context of the standalone PNP equations, numerous structure-preserving schemes have been proposed. Researchers have developed finite difference, finite element, and finite volume methods that successfully maintain positivity and energy stability, often relying on exponential transformations, Slotboom approximations, or sophisticated stabilized implicit treatments \cite{flavell2014conservative, gao2017linearized, gao2018linearized, hu2020fully, liu2021positivity, prohl2009convergent, shen2021unconditionally}.

Extending these structure-preserving techniques to the fully coupled PNPNS system introduces additional complexities due to the convective coupling and the electric body force. Early computational approaches often rely on fully coupled, implicit discretizations to guarantee energy stability, which inherently require the solution of massive, highly nonlinear algebraic systems at each time step, resulting in prohibitive computational costs \cite{bauer2012stabilized, liu2017efficient, metti2016energetically, prohl2010convergent,linga2020transient}. To mitigate this computational bottleneck, a recent trend in the literature has been the development of decoupled, linearized, and energy-stable schemes. Notably, the scalar auxiliary variable (SAV) approach and its variants have been widely adopted to explicitly treat the nonlinear coupling terms while artificially preserving a modified energy dissipation law \cite{shen2019new, zhou2023efficient}. Building upon these decoupling strategies, several efficient time-stepping schemes, including pressure-correction and projection methods, have been recently advanced for the electro-hydrodynamic system, achieving high-order temporal accuracy and unconditional energy stability \cite{dehghan2023optimal, he2025stability, pan2024linear, yu2025decoupled}.

Despite these significant advancements, dealing with the incompressibility constraint of the fluid in a structurally consistent manner remains a critical issue. Standard finite element methods often fail to yield strictly divergence-free velocity fields, which can lead to the violation of exact mass conservation in the advection of ionic concentrations and induce spurious oscillations. To address this, fractional-step projection methods have been thoroughly investigated for incompressible flows \cite{guermond2006overview}, and the rotational form of the pressure-correction scheme has proved particularly effective in minimizing artificial boundary layers \cite{timmermans1996approximate}. Moreover, to strictly satisfy the incompressibility condition at the discrete level, mixed finite element methods and virtual element methods have been introduced into the discretization of both the Stokes/Navier-Stokes equations and the coupled electro-hydrodynamic systems \cite{dehghan2023optimal,correa2024banach, correa2023mixed, he2018mixed, he2021mixed}. However, literature proposing a fully decoupled, structure-preserving mixed finite element method for the complete PNPNS system, backed by fully discrete structural proofs, remains scarce.

This paper proposes a novel decoupled, structure-preserving discrete method for PNPNS equations in the mixed finite element framework. In our non-dimensionalized formulation, the Navier-Stokes equations are cast in a rotational form that naturally emphasizes vorticity dynamics and enables the straightforward construction of structure-preserving numerical discretizations. From a scheme construction perspective, it is a fully linearized and decoupled numerical scheme that allows the PNP and NS subsystems to be solved separately as linear systems. An additional key advantage is that it naturally leads to a leap-frog-inspired time integration strategy, which handles all nonlinearities in the PNPNS system by staggering the time evolution of fluid and electrochemical variables. As a result, the fully discrete scheme generates only linear, decoupled algebraic systems that can be solved independently. We also rigorously prove that the proposed method strictly preserves mass, positive and negative ions, and guarantees the energy stability at the fully discrete level. Comprehensive numerical experiments are conducted to validate our theoretical claims. 

The remainder of this paper is organized as follows. Section~\ref{Sec: scheme} starts with the PNPNS model and its properties. It then introduces the functional setting and also details the construction of the decoupled, structure-preserving mixed finite element scheme. Next, in Section~\ref{Sec proof}, we provide the rigorous theoretical proofs for mass conservation and the energy stability of the proposed scheme. Section~\ref{SEC: numeric} presents extensive numerical experiments to validate the theoretical findings and demonstrate the accuracy and robustness. 

\section{Model and discretization}
\label{Sec: scheme}
\subsection{The PNPNS model and its reformulation}
Let \(\Omega \subset \mathbb{R}^d\), $d\in\left\lbrace2,3\right\rbrace$, be a bounded contractible domain with a Lipschitz boundary \(\partial \Omega\). 
In a space-time domain $\Omega\times(0, T)$, we consider the following normalized nondimensional coupled system of PNP equations and NS equations, \cite{yu2025decoupled},
\begin{subequations} \label{Eq:PNP-NS}
\begin{align}
	\frac{\partial p}{\partial t} + \left(\boldsymbol{u} \cdot \nabla\right)p &= \nabla \cdot \left (\nabla p + p\nabla \psi \right), \label{Eq:PNP-NSa} \\
	\frac{\partial n}{\partial t} + \left(\boldsymbol{u} \cdot \nabla \right)n &= \nabla \cdot \left(\nabla n - n\nabla \psi\right), \label{Eq:PNP-NSb} \\
	-\varepsilon \Delta \psi &= p - n, \label{Eq:PNP-NSc} \\
	\frac{\partial\boldsymbol{u}}{\partial t} + \nabla \times \boldsymbol{\omega} + \nabla P &= -\nabla \psi\left(p - n \right), \label{Eq:PNP-NSd} \\
	\boldsymbol{\omega} &= \nabla \times \boldsymbol{u}, \label{Eq:PNP-NSe} \\
	\nabla \cdot \boldsymbol{u} &= 0, \label{Eq:PNP-NSf}
\end{align}
\end{subequations}
where $\boldsymbol{u} $ is the velocity field, $\boldsymbol{\omega}$ is the vorticity field,  the total pressure P is given by $p_{\text{static}} + \tfrac{1}{2} \boldsymbol{u}^2$ where $p_{\text{static}}$ denotes the static pressure. Besides, the variables \(p\) and \(n\) represent the concentration functions of positive and negative ions in the fluid. \(\psi\) is the electric potential. And the parameter $\varepsilon$ is the dielectric permittivity. Without loss of generality, we will use $\varepsilon =1$ in the remainder of the paper.  

To reveal the underlying thermodynamic structure and to facilitate the development of the energy-stable numerical scheme, we reformulate the above system by introducing chemical potentials for the macroscopic ionic species,
\[\mu = \ln p + \psi, \quad \nu = \ln n - \psi.\]
By direct differentiation, the coupled diffusion and electromigration fluxes can be compactly rewritten as $\nabla p + p\nabla \psi = p\nabla \mu$ and $\nabla n - n\nabla \psi = n\nabla \nu$.
Furthermore, this transformation allows us to handle the singular electrical coupling term $-\nabla \psi \left(p - n\right)$ as follows. Summing the potential gradients leads to 
\[p\nabla \mu + n\nabla \nu = \nabla \left(p + n \right) + \nabla \psi\left(p - n \right),\] 
which implies 
\[-\nabla \psi \left(p - n \right) = -\left(p\nabla \mu + n\nabla \nu \right) + \nabla \left(p + n \right).\]
Substituting this relation into the Navier-Stokes equations and introducing a modified effective pressure \begin{equation} \label{Eq: phi}
\phi = P - p - n
\end{equation}
to absorb the pure gradient term $\nabla \left(p + n \right)$ via the incompressibility constraint, we arrive at the following reformulated PNPNS system,
\begin{subequations} \label{Eq: PNPNS}
\begin{align}
	\frac{\partial p}{\partial t} + \left(\boldsymbol{u} \cdot \nabla\right)p - \nabla \cdot \left(p \nabla \mu \right)&=0,\label{Eq: PNPNSa}\\
	\frac{\partial n}{\partial t} + \left(\boldsymbol{u} \cdot \nabla \right)n - \nabla \cdot \left (n \nabla \nu \right)&=0,\label{Eq: PNPNSb}\\
	- \nabla \cdot \nabla \psi &= p - n,\label{Eq: PNPNSc}\\
	\frac{\partial \boldsymbol{u}}{\partial t} +\boldsymbol{\omega} \times \boldsymbol{u} + \nabla \times \boldsymbol{\omega} + \nabla \phi &= -\left(p \nabla \mu + n \nabla \nu \right),\label{Eq: PNPNSd}\\
	\boldsymbol{\omega} &= \nabla \times \boldsymbol{u},\label{Eq: PNPNSe}\\
	\nabla \cdot \boldsymbol{u} &= 0.\label{Eq: PNPNSf}
\end{align}
\end{subequations}
This reformulated system maintains strict physical equivalence with the original equations, but explicitly isolates the thermodynamically consistent driving forces $\nabla \mu$ and $\nabla \nu$, and thus is highly advantageous for establishing the energy dissipation law.
The system is supplemented with one of the two sets of boundary conditions,
\begin{equation}\label{Eq: boundary conditions}
\boldsymbol{u} \cdot \boldsymbol{n} = 0,\quad \boldsymbol{u} \times \boldsymbol{n} = 0,\quad \left(p \nabla \mu \right)  \cdot \boldsymbol{n} = 0,\quad \left (n \nabla \nu \right) \cdot \boldsymbol{n} = 0, \quad \nabla \psi \cdot \boldsymbol{n} = 0,  \quad\text{on}\, \partial \Omega \times (0, T],
\end{equation}
or periodic boundary conditions. Additionally, The system is closed subject to initial conditions, 
\begin{equation}\label{Eq: initial conditions}
p(\boldsymbol{x}, 0) = p_0(\boldsymbol{x}),\quad n(\boldsymbol{x}, 0) = n_0(\boldsymbol{x}),\quad \boldsymbol{u}(\boldsymbol{x}, 0) = \boldsymbol{u}_0(\boldsymbol{x}).
\end{equation}
Note that, one can compute the initial condition $\psi_0$ by solving the Poisson problem \eqref{Eq: PNPNSc} using the boundary conditions and initial conditions. The initial vorticity $\boldsymbol{\omega}_0$ can be obtained from \eqref{Eq: PNPNSe} with the initial velocity $\boldsymbol{u}_0$.

For the given boundary condition $\boldsymbol{u} \cdot \boldsymbol{n} = 0 $  or periodic boundary conditions, the conservation of ion mass holds, i.e., for any time \( t \in (0, T] \), \cite{yu2025decoupled},
\[
\int_{\Omega} p(x,t) \, \mathrm{d}x = \int_{\Omega} p(x,0) \, \mathrm{d}x, \qquad \int_{\Omega} n(x,t) \, \mathrm{d}x = \int_{\Omega} n(x,0) \, \mathrm{d}x.
\]
Moreover, the PNPNS system possesses a key energy dissipation property of the form
\[
\frac{\mathrm{d}\mathcal{E}}{\mathrm{d}t}  = -\int_{\Omega} \boldsymbol{\omega}^2 + p\left(\nabla \mu\right)^2 + n\left(\nabla v\right)^2  \, \mathrm{d}x < 0, 
\]
where the total energy functional $ \mathcal{E}$  is defined  as
\[
\mathcal{E} = \int_{\Omega} p\left(\ln p - 1\right) + n\left(\ln n - 1
\right) + \frac{1}{2}\left(\nabla \psi\right)^2 + \frac{1}{2}\boldsymbol{u}^2  \, \mathrm{d}x.
\] 
\subsection{Notations and the weak formulation}

The symbol \(\langle \cdot, \cdot\rangle_\Omega\) denotes the standard inner product on \(L^2(\Omega)\), and \(\|\cdot\| \) denotes the $L^2$ norm. And we will take three dimensions ($d=3$) for explaining the present method. And finding its two-dimensional version is straightforward.

Recall the following well-known functional spaces,
\[
\begin{aligned}
&L^2(\Omega) := \left\{ f \mid \langle f, f \rangle_\Omega < +\infty \right\},\\
&\boldsymbol{H}(\operatorname{curl};\Omega) := \left\{ \boldsymbol{\omega} \in L^2(\Omega) \mid \nabla \times \boldsymbol{\omega} \in [L^2(\Omega)]^3 \right\}, \\
&\boldsymbol{H}(\operatorname{div};\Omega) := \left\{ \boldsymbol{u} \in [L^2(\Omega)]^3 \mid \nabla \cdot \boldsymbol{u} \in L^2(\Omega) \right\}, \\
&H^1(\Omega) := \left\{ \phi \in L^2(\Omega) \mid \nabla \phi \in [L^2(\Omega)]^2 \right\}.
\end{aligned}
\]
They form a well-known Hilbert complex, i.e., the de Rham complex \cite{arnold2010finite,arnold2006finite, bochev2003discourse, palha2017mass}, 
\begin{equation}
\mathbb{R} \hookrightarrow H^1(\Omega) \xrightarrow{\nabla } \boldsymbol{H}(\mathrm{curl};\Omega) \xrightarrow{\nabla \times} \boldsymbol{H}(\mathrm{div};\Omega) \xrightarrow{\nabla \cdot} L^2(\Omega) \to 0.
\end{equation}
We denote by $\mathcal{T}$ the trace operator, which maps a function to its restriction on a boundary segment $\Gamma \subseteq \partial\Omega$. 
We then introduce the trace spaces as
\[
\mathcal{T} \boldsymbol{H}(\text{curl}; \Omega, \Gamma) := \left\{ \mathcal{T} \boldsymbol{\omega} \mid \boldsymbol{\omega} \in \boldsymbol{H}(\text{curl}; \Omega) \right\},
\]
\[
\mathcal{T} H(\text{div}; \Omega, \Gamma) := \left\{ \mathcal{T} \boldsymbol{u} \mid \boldsymbol{u} \in \boldsymbol{H}(\text{div}; \Omega) \right\}.
\]
\[
H^{1/2}(\Omega, \Gamma) := \left\{ \mathcal{T} \phi \mid \phi \in H^1(\Omega) \right\}.
\]
A complete introduction to Sobolev spaces can be found in \cite{oden2017applied}.

We employ the notations  $G(\Omega)$, $\boldsymbol{C}(\Omega)$, $\boldsymbol{D}(\Omega)$, and $S(\Omega)$ to express finite-dimensional function spaces which are subsets of Sobolev spaces, namely,
\[G(\Omega) \subset H^1(\Omega),\ \boldsymbol{C}(\Omega) \subset H(\text{curl};\Omega),\ \boldsymbol{D}(\Omega) \subset H(\text{div};\Omega),\  S(\Omega) \subset L^2(\Omega).\\
\]
And form discrete de Rham complexes in three-dimensional space,
\begin{equation} \label{Eq: de}
\mathbb{R} \hookrightarrow G(\Omega) \xrightarrow{\nabla} \boldsymbol{C}(\Omega) \xrightarrow{\nabla  \times} \boldsymbol{D}(\Omega) \xrightarrow{\nabla \cdot} S(\Omega) \to 0.
\end{equation}
An additional required regularity is that these discrete spaces exhibit sufficient regularity to ensure the $L^2$-integrability for the nonlinear terms at the semi-discrete formulation \eqref{eq: formula}. This usually is the case for particular discrete spaces such as the ones we will use for numerical experiments in Section~\ref{SEC: numeric}, namely, the mimetic spectral element spaces. 
We also introduce the subspace, the zero-trace subspace of 
\( \boldsymbol{D}(\Omega) \), comprises the element in \( \boldsymbol{D}(\Omega) \) whose trace vanishes on $\partial\Omega$, i.e.,
\[
\boldsymbol{D} _0(
\Omega, \Gamma) := \left\{ \boldsymbol{u}_h \mid \boldsymbol{u}_h \in \boldsymbol{D}(\Omega),\ \mathcal{T}\boldsymbol{u}_h = 0 \text{ on } \partial\Omega\right\}
.
\]
Throughout the paper, we use subscript \( h \) to denote elements of finite-dimensional spaces. 

Upon these notations, we now propose a semi-discrete formulation of \eqref{Eq: PNPNS}: Seek
$(p_h, n_h, \psi_h, \boldsymbol{u_h}, \boldsymbol{\omega_h}, \phi_h) \in \left[G(\Omega) \right]^3\times \boldsymbol{D}_0(\Omega) \times \boldsymbol{C}( \Omega) \times S(\Omega),$
such that,$\forall (q_h, m_h, \varphi_h, \boldsymbol{v_h},\boldsymbol{w_h}, \alpha_h) \in \left[G(\Omega) \right]^3 \times \boldsymbol{D}_0(\Omega) \times \boldsymbol{C}(\Omega) \times S(\Omega),$
\begin{subequations}\label{eq: formula}
\begin{align}
	\left\langle\frac{\partial p_h}{\partial t},q_h\right\rangle+\left\langle p_h \nabla q_h,\boldsymbol{u}_h\right\rangle +\left\langle p_h\nabla\mu_h,\nabla q_h\right\rangle&=0, \\
	\left\langle\frac{\partial n_h}{\partial t},m_h\right\rangle+\left\langle n_h\nabla m_h,\boldsymbol{v}_h\right\rangle  +\left\langle n_h\nabla\nu_h,\nabla m_h\right\rangle&=0, \\	
	\left\langle\nabla\psi_h,\nabla\varphi_h\right\rangle & =\left\langle p_h-n_h,\varphi_h\right\rangle, \\
	\left\langle \mu_h,\xi_h\right\rangle&=	\left\langle\ln p_h+\psi_h,\xi_h\right\rangle\\
	\left\langle \nu_h,\eta_h\right\rangle&=	\left\langle\ln n_h-\psi_h,\eta_h\right\rangle\\
	\left\langle\frac{\partial\boldsymbol{u}_h}{\partial t},\boldsymbol{v}_h\right\rangle+\left\langle\boldsymbol{\omega}_h\times\boldsymbol{u}_h,\boldsymbol{v}_h\right\rangle+\left\langle\nabla\times\boldsymbol{\omega}_h,\boldsymbol{v}_h\right\rangle-\left\langle\phi_h,\nabla\cdot\boldsymbol{v}_h\right\rangle \hspace{-2cm}& \label{formula f}\\&=-\left\langle p_h\nabla\mu_h,\boldsymbol{v}_h\right\rangle-\left\langle n_h\nabla\nu_h,\boldsymbol{v}_h\right\rangle,\notag \\
	\left\langle\boldsymbol{\omega}_h,\boldsymbol{w}_h\right\rangle-\left\langle\boldsymbol{u}_h,\nabla\times\boldsymbol{w}_h\right\rangle & =0, \label{formula g}\\
	\langle\nabla\cdot\boldsymbol{u}_h,\alpha_h\rangle &=0,\label{formula h}
\end{align}
\end{subequations}
subject to initial conditions in \eqref{Eq: initial conditions}.
Note that the essential boundary condition $\boldsymbol u_h\cdot\boldsymbol{n_h}=0$ is enforced by selecting $\boldsymbol{u}_h\in \boldsymbol{D}_0(\mathrm{div};\Omega)$, and other boundary conditions are included in a natural way through boundary integrals that have vanished.

\subsection{A decoupled numerical scheme}
\label{subsec1}
To decouple and to linearize the semi-discrete formulation \eqref{eq: formula}, we propose the following temporal scheme.

Suppose we have generated a temporal sequence 
$$
t^0,\ t^{\frac{1}{2}},\ t^1,\ t^{\frac{3}{2}},\ \cdots,\ t^{k-\frac{1}{2}},\ t^k,\ t^{k+\frac{1}{2}},\ \cdots,
$$
where $t^0=0$, $t^k-t^{k-1}=\varDelta t$ $\left(k\in\left\lbrace1,2,3,\cdots\right\rbrace\right)$ is a constant, and $t^{k-\frac{1}{2}} = \dfrac{t^{k-1}+t^k}{2}$. And we will use the superscript  $k$ to denote a temporally discrete variable at the time instant $t^k$. For example, $p_h^k=p_h(\boldsymbol{x},t^k)$. 

Using the initial conditions, we can find the approximations $p_h^0$, $n_h^0 $, and $\boldsymbol{u}_h^{0}$. We then can compute other initial conditions, $\mu_h^0$, $\nu_h^0$, $\psi_h^0$, $\boldsymbol{\omega}_h^0$, through
\begin{subequations}\label{eq sche 1}
\begin{align}
	\left\langle\mu_h^0,\xi_h\right\rangle&=\left\langle \ln p_h^{0} + \psi_h^{0},\xi_h\right\rangle,\  \forall \xi_h \in G(\Omega),\label{eq sche a}\\
	\left\langle\nu_h^0,\eta_h\right\rangle&=\left\langle \ln n_h^{0} - \psi_h^{0},\eta_h\right\rangle,\  \forall \eta_h \in G(\Omega),\label{eq sche b}\\
	\left\langle \nabla \psi_h^0, \nabla \varphi_h \right\rangle &= \left\langle p_h^0 - n_h^0, \varphi_h \right\rangle,\  \forall \varphi_h \in G(\Omega),\label{eq sche c}\\
	\left\langle \boldsymbol{\omega}_h^{0}, \boldsymbol{w_h} \right\rangle &= \left\langle \boldsymbol{u}_h^{0}, \nabla \times \boldsymbol{w_h} \right\rangle,\ \forall \boldsymbol{w_h} \in \boldsymbol{C}(\Omega)\label{eq sche d}
\end{align}
\end{subequations}
To start the regular iterations, we will first need to find $\boldsymbol{u}_h^{\frac{1}{2}}$ and $\boldsymbol{\omega}_h^{\frac{1}{2}}$ by solving the Navier-Stokes part of \eqref{eq: formula}, i.e., \eqref{formula f} - \eqref{formula h}, using, for example, an explicit temporal scheme for the time step from $t^0$ to $t^\frac{1}{2}$.

Then, the regular iterations of the fully discrete scheme are as follows: For $k=1,2,3,\cdots$, sequentially, \\
\textbf{(step 1)} seek $\left(p_h^k, \ n_h^k,\ \psi_h^k \right) \in \left[G(\Omega)\right]^3$, such that,
$\forall \left(q_h, m_h, \varphi_h\right)\in \left[G(\Omega)\right]^3 $,
\begin{subequations}\label{eq scheme 1}
\begin{align}
	& \left\langle \frac{p_h^k - p_h^{k-1}}{\varDelta t}, q_h \right\rangle- \left\langle \frac{p_h^{k-1} + p_h^k}{2} \nabla q_h, \boldsymbol{u}_h^{k-\frac{1}{2}} \right\rangle +\left\langle \frac{p_h^{k-1} + p_h^k}{2} \nabla \mu_h^{k-1}, \nabla q_h \right\rangle=0, \label{eq scheme 1a}\\
	& \left\langle \frac{n_h^k - n_h^{k-1}}{\varDelta t}, m_h \right\rangle - \left\langle \frac{n_h^{k-1} + n_h^k}{2} \nabla m_h, \boldsymbol{u}_h^{k-\frac{1}{2}} \right\rangle +\left\langle \frac{n_h^{k-1} + n_h^k}{2} \nabla \nu_h^{k-1}, \nabla m_h \right\rangle=0,\label{eq scheme 1b} \\
	& \left\langle \nabla \psi_h^k, \nabla \varphi_h \right\rangle = \left\langle p_h^k - n_h^k, \varphi_h \right\rangle,\label{eq scheme 1c}
\end{align}
\end{subequations}
\textbf{(step 2)} seek $\left(\mu_h^k,\ \nu_h^k \right) \in \left[G(\Omega)\right]^2$, such that, $\forall \left(\xi_h,\eta_h\right)\in \left[G(\Omega)\right]^2$,
\begin{subequations}\label{eq scheme 2}
\begin{align}
	&\left\langle\mu_h^{k},\xi_h\right\rangle = \left\langle\ln p_h^{k} + \psi_h^{k},\xi_h \right\rangle, \label{eq scheme 2a}\\
	&\left\langle\nu_h^{k},\eta_h\right\rangle = \left\langle\ln n_h^{k} - \psi_h^{k},\eta_h \right\rangle,\label{eq scheme 2b}
\end{align}
\end{subequations}	
\textbf{(step 3)} seek $ \left(\boldsymbol{u}_h^{k+\frac{1}{2}},\ \boldsymbol{\omega}_h^{k+\frac{1}{2}},\phi_h^1\right)\in \boldsymbol{D}_0(\Omega) \times \boldsymbol{C}(\Omega)\times S(\Omega)$, such that, $\forall \left(\boldsymbol{v}_h,\boldsymbol{w_h}, \alpha_h\right)\in \boldsymbol{D}_0(\Omega) \times \boldsymbol{C}(\Omega)\times S(\Omega)$,
\begin{subequations}\label{eq scheme 3}
\begin{align}
	&\left\langle \frac{\boldsymbol{u}_h^{k+\frac{1}{2}} - \boldsymbol{u}_h^{k-\frac{1}{2}}}{\varDelta t}, \boldsymbol{v}_h \right\rangle + \left\langle \boldsymbol{\omega}_h^{k-\frac{1}{2}} \times \boldsymbol{u}_h^{k+\frac{1}{2}}, \boldsymbol{v}_h \right\rangle + \left\langle \nabla \times \boldsymbol{\omega}_h^{k+\frac{1}{2}}, \boldsymbol{v}_h \right\rangle - \left\langle \phi_h^k, \nabla \cdot \boldsymbol{v}_h \right\rangle \label{eq scheme 3a} 
	\\ 
	&\hspace{3cm}= -\left\langle \frac{p_h^{k-1} + p_h^k}{2} \nabla \mu_h^k, \boldsymbol{v}_h \right\rangle - \left\langle \frac{n_h^{k-1} + n_h^k}{2} \nabla \nu_h^k, \boldsymbol{v}_h \right\rangle, \notag
	\\
	&\left\langle \boldsymbol{\omega}_h^{k+\frac{1}{2}}, \boldsymbol{w_h} \right\rangle - \left\langle \boldsymbol{u}_h^{k+\frac{1}{2}}, \nabla \times \boldsymbol{w_h} \right\rangle = 0, \label{eq scheme 3b} 
	\\
	&\left\langle \nabla \cdot \boldsymbol{u}_h^{k+\frac{1}{2}}, \alpha_h \right\rangle = 0,  \label{eq scheme 3c}
\end{align}
\end{subequations}
until $t^{k+\frac{1}{2}} \geq T$ or any other criterion is reached. 

\section{Properties of the numerical scheme}
\label{Sec proof}

We show below that our decoupled numerical scheme, \eqref{eq scheme 1} - \eqref{eq scheme 3}, enjoys properties, i.e.
mass conservation and conditional energy stability. Note that we use \(C\) to denote the independent constant, and it may take different values in different occurrences.
\begin{theorem}\label{t1}\emph{
	For all \(k\in\left\lbrace1,2,3,\cdots\right\rbrace\), given \( \left(p_{h}^{k-1},\ n_{h}^{k-1},\ \boldsymbol{u}_{h}^{k-\frac{1}{2}},\ \phi_{h}^{k-1} \right) \in \left[G(\Omega)\right]^2\times \boldsymbol{D}_0(\Omega) \times S(\Omega) \), if the solution $p_{h}^{k}$, $n_{h}^{k}$ in \( \Omega \) is positive and bounded, i.e., there exist positive constants $m$ and $M $  independent of the mesh size $h$ and time step index $k$, such that \(0<m < p_{h}^{k}\leq M\), and \(0<m <n_{h}^{k}<M\),  then the proposed scheme (\ref{eq scheme 1})-(\ref{eq scheme 3}) has the following properties.
	\begin{itemize}
		\item \textbf{Mass Conservation}
		\[
		\left\langle p_{h}^{k}, 1 \right\rangle = \left\langle p_{h}^{k-1}, 1 \right\rangle, \quad \left\langle n_{h}^{k}, 1 \right\rangle = \left\langle n_{h}^{k-1}, 1 \right\rangle.
		\]
		\item \textbf{Conditional Energy Stability}\\
		The discrete energy at step \(k\) is defined as:
		\[
		\mathcal{E}_{h}^{k} = \left\langle p_h^k, \ln p_h^{k}-1 \right\rangle + \left\langle n_h^k, \ln n_h^{k}-1 \right\rangle + \frac{1}{2} \left\|\nabla \psi_h^k\right\|^2 + \frac{1}{2} \left\|\boldsymbol{u}_h^{k+\frac{1}{2}} \right\|^2.
		\]
		There exists \(\varDelta t = \min\left\{ \frac{m}{2M^2}, {C h^2} \right\}\) such that: 
		\begin{equation*}
			\begin{aligned}
				\mathcal{E}_{h}^{k}-\mathcal{E}_{h}^{k-1} &=\left\langle  p_h^{k}\left(\ln p_h^k-1\right) - p_h^{k-1}\left(\ln p_h^{k-1} - 1\right),1 \right\rangle +\left\langle  n_h^{k}(\ln n_h^k-1) - n_h^{k-1}\left(\ln n_h^{k-1} - 1\right),1 \right\rangle\\&\quad+ \frac{1}{2} \left(\left \| \nabla \psi_h^k \right\|^2 - \left\| \nabla \psi_h^{k-1}\right \|^2 \right)+\frac{1}{2} \left( \left\| \boldsymbol{u}_h^{k+\frac{1}{2}} \right\|^2 - \left\| \boldsymbol{u}_h^{k-\frac{1}{2}} \right\|^2 \right)\\
				&\leq \left( -\frac{1}{2} +  C \varDelta t h^{-2}\right) \left\|\nabla \psi_h^k - \nabla \psi_h^{k-1}\right\|^2
				- \varDelta t\left\|\boldsymbol{\omega}_h^{k+\frac{1}{2}}\right\|^2  - \frac{\varDelta t}{2} M \left\Vert  \nabla \mu_h^{k-1} \right\Vert^2- \frac{\varDelta t}{2} M \left\Vert  \nabla \nu_h^{k-1} \right\Vert^2  \\
				&\quad+\varDelta t \left( \varDelta t M^2 - \frac{m}{2} \right) \left\|\nabla \mu_h^k\right\|^2  + \varDelta t \left( \varDelta t M^2 - \frac{m}{2} \right) \left\|\nabla \nu_h^k\right\|^2  \\
				&\quad+ \left(C  \varDelta t h^{-2}-C'\right) \left\|p_h^k - p_h^{k-1}\right\|^2 + \left(C \varDelta t h^{-2}-C'\right)\left\|n_h^k - n_h^{k-1}\right\|^2
				\\&\leq0.
			\end{aligned}
		\end{equation*}
	\end{itemize}
}\end{theorem}
\paragraph{Proof} 
Mass Conservation follows directly by choosing test functions $q_h = 1,\ m_h = 1$ in \eqref{eq scheme 1a} and \eqref{eq scheme 1b}, respectively.
To prove the conditional energy stability,
we first express
\[
\begin{aligned}
\mathcal{E}_{h}^{k}-\mathcal{E}_{h}^{k-1} = E_1 + E_2 + E_3,
\end{aligned}
\]
where \[
\begin {aligned}
& E_1=\left[\left \langle  p_h^{k}\left(\ln p_h^k-1\right) - p_h^{k-1}\left(\ln p_h^{k-1} - 1\right),1 \right\rangle +\left\langle  n_h^{k}\left(\ln n_h^k-1\right) - n_h^{k-1}\left(\ln n_h^{k-1} - 1\right),1 \right\rangle \right],\\
& E_2= \frac{1}{2} \left( \left\| \nabla \psi_h^k \right\|^2 - \left\| \nabla \psi_h^{k-1} \right\|^2 \right), \\
& E_3 =\frac{1}{2} \left( \left\| \boldsymbol{u}_h^{k+\frac{1}{2}} \right\|^2 - \left\| \boldsymbol{u}_h^{k-\frac{1}{2}} \right\|^2 \right).
\end {aligned}
\]
To estimate \(E_1\), from the properties of function $x\left(\ln x-1\right)$, we know there must exist a positive constant $C'$ such that

\[
\begin{aligned}
\left\langle p_h^k - p_h^{k-1}, \ln p_h^{k} \right\rangle- C' \left\|p_h^k - p_h^{k-1}\right\|^2 &\geq \left\langle  p_h^{k}\left(\ln p_h^k-1\right) - p_h^{k-1}\left(\ln p_h^{k-1} - 1\right),1 \right\rangle ,\\
\left\langle n_h^k - n_h^{k-1}, \ln n_h^{k} \right\rangle- C' \left\|n_h^k - n_h^{k-1}\right\|^2 &\geq \left\langle  n_h^{k}\left(\ln n_h^k-1\right) - n_h^{k-1}\left(\ln n_h^{k-1} - 1\right),1 \right\rangle.
\end{aligned}
\]
See \ref{App: A} for more details on these inequalities. 
Next, to estimate $\left\langle p_h^k - p_h^{k-1}, \ln p_h^{k} \right\rangle$  and $\left\langle n_h^k - n_h^{k-1} ,\ln n_h^{k} \right\rangle$, since \eqref{eq scheme 1a} and \eqref{eq scheme 1b} hold for all test functions  in $ G(\Omega)$, we can select the test functions to  be \( q_h = \mu_h^{k} \)  and \( m_h = \nu_h^{k} \), respectively. As a result, we obtain
\[
\begin{aligned}
\left\langle \frac{p_h^k - p_h^{k-1}}{\varDelta t}, \ln p_h^k + \psi_h^k \right\rangle - \left\langle \frac{p_h^{k-1} + p_h^k}{2} \nabla \mu_h^{k}, \boldsymbol{u}_h^{k-\frac{1}{2}} \right\rangle &=-  \left\langle \frac{p_h^{k-1} + p_h^{k}}{2}\nabla \mu_h^{k-1},  \nabla \mu_h^{k},  \right\rangle,
\\
\left\langle \frac{n_h^k - n_h^{k-1}}{\varDelta t}, \ln n_h^k - \psi_h^k \right\rangle - \left\langle \frac{n_h^{k-1} + n_h^k}{2} \nabla \nu_h^{k}, \boldsymbol{u}_h^{k-\frac{1}{2}} \right\rangle &=- \left\langle \frac{n_h^{k-1} + n_h^{k}}{2}\nu_h^{k-1},  \nabla \nu_h^{k} \right\rangle.
\end{aligned}
\]
They further lead to
\[
\begin{aligned}
\left\langle p_h^k - p_h^{k-1}, \ln p_h^k\right\rangle &=\varDelta t \left\langle \frac{p_h^{k-1} + p_h^{k}}{2} \nabla \mu_h^k,  \boldsymbol{u}_h^{k-\frac{1}{2}} \right\rangle - \varDelta t \left\langle \frac{p_h^{k-1} + p_h^k}{2} \nabla \mu^{k-1}, \nabla \mu_h^k \right\rangle-\left\langle p_h^k - p_h^{k-1}, \psi_h^k \right\rangle,
\\
\left\langle n_h^k - n_h^{k-1}, \ln n_h^k \right \rangle &=  \varDelta t \left\langle \frac{n_h^{k-1} + n_h^{k}}{2} \nabla \nu_h^k,  \boldsymbol{u}_h^{k-\frac{1}{2}} \right\rangle - \varDelta t \left\langle \frac{n_h^{k-1} + n_h^k}{2} \nabla \nu^{k-1}, \nabla \nu_h^k \right\rangle+\left\langle n_h^k - n_h^{k-1}, \psi_h^k \right\rangle.
\end{aligned}
\]
With these relations, we can find
\begin{equation}\label{E1}
\begin{aligned}
	E_1&=\left\langle  p_h^{k}\left(\ln p_h^k-1\right) - p_h^{k-1}\left(\ln p_h^{k-1} - 1\right),1 \right\rangle +\left\langle  n_h^{k}\left(\ln n_h^k-1\right) - n_h^{k-1}\left(\ln n_h^{k-1} - 1\right),1 \right\rangle  \\&\leq  \varDelta t \left\langle \frac{p_h^{k-1} + p_h^{k}}{2} \nabla \mu_h^k,  \boldsymbol{u}_h^{k-\frac{1}{2}} \right\rangle - \varDelta t \left\langle \frac{p_h^{k-1} + p_h^k}{2} \nabla \mu_h^{k-1}, \nabla \mu_h^k \right\rangle \\&\qquad+\varDelta t \left\langle \frac{n_h^{k-1} + n_h^{k}}{2} \nabla \nu_h^k,  \boldsymbol{u}_h^{k-\frac{1}{2}} \right\rangle - \varDelta t \left\langle \frac{n_h^{k-1} + n_h^k}{2} \nabla \nu_h^{k-1}, \nabla \nu_h^k \right\rangle \\&\qquad-\left\langle p_h^k - p_h^{k-1}, \psi_h^k \right\rangle+\left\langle n_h^k - n_h^{k-1}, \psi_h^k \right\rangle - C' \left\|p_h^k - p_h^{k-1}\right\|^2 - C' \left\|n_h^k - n_h^{k-1}\right\|^2,
\end{aligned}
\end{equation}
where $-\left\langle p_h^k - p_h^{k-1}, \psi_h^k \right\rangle+\left\langle n_h^k - n_h^{k-1}, \psi_h^k\right \rangle =	- \left\langle \left(p_h^k - n_h^k\right) - \left(p_h^{k-1} - n_h^{k-1}\right), \psi_h^{k}\right \rangle$,
which allows us to bound the sum of the two left-hand terms by equivalently estimating the right-hand term. Next, we select the test function in \eqref{eq scheme 1c} to be \(\varphi_h = \psi_h^k \) and get
\[
\left\langle p_h^k - n_h^k,\psi_h^k \right\rangle= \left\| \nabla \psi_h^k\right \|^2 .
\]
We can do the same to the previous time step. Thus we also have
\[
\left\langle p_h^{k-1} - n_h^{k-1}, \psi_h^k \right\rangle=\langle \nabla \psi_h^{k-1},  \nabla \psi_h^k\rangle.
\]
With the inequality \(-a(a - b) = -\frac{1}{2}\left(a^2 - b^2 + (a - b)^2\right)\) and we can conclude that
\begin{equation}\label{E2}
\begin{aligned}
	-\left\langle p_h^k - p_h^{k-1}, \psi_h^k \right\rangle+\left\langle n_h^k - n_h^{k-1}, \psi_h^k \right\rangle &=-\left\| \nabla \psi_h^k \right\|^2+\left\langle \nabla \psi_h^{k-1},  \psi_h^k\right\rangle\\&=-\frac{1}{2} \left( \left\|\nabla \psi_h^{k}\right\|^2 - \left\|\nabla \psi_h^{k-1}\right\|^2 + \left\|\nabla (\psi_h^{k} - \psi_h^{k-1})\right\|^2 \right).
\end{aligned}
\end{equation}
Substituting \eqref{E2} into \eqref{E1} and using the definition of $E_2$, we get
\begin{equation}\label{eq E3}
\begin{aligned}
	E_1+E_2&\leq  \varDelta t \left\langle \frac{p_h^{k-1} + p_h^{k}}{2} \nabla \mu_h^k,  \boldsymbol{u}_h^{k-\frac{1}{2}} \right\rangle - \varDelta t \left\langle \frac{p_h^{k-1} + p_h^k}{2} \nabla \mu_h^{k-1}, \nabla \mu_h^k \right\rangle \\&\qquad+\varDelta t \left\langle \frac{n_h^{k-1} + n_h^{k}}{2} \nabla \nu_h^k,  \boldsymbol{u}_h^{k-\frac{1}{2}} \right\rangle - \varDelta t \left\langle \frac{n_h^{k-1} + n_h^k}{2} \nabla \nu_h^{k-1}, \nabla \nu_h^k \right\rangle \\& \qquad- \frac{1}{2}\left\|\nabla (\psi_h^{k} - \psi_h^{k-1})\right\|^2 - C' \left\|p_h^k - p_h^{k-1}\right\|^2 - C' \left\|n_h^k - n_h^{k-1}\right\|^2.	
\end{aligned}
\end{equation}
To evaluate the term $E_3$, we now select the test function in \eqref{eq scheme 3a} to be $\boldsymbol{v}_h = \boldsymbol{u}_h^{k+\frac{1}{2}}$ since \eqref{eq scheme 3a} holds $\forall \boldsymbol{v}_h \in \boldsymbol{D}(\Omega)$. This gives
\begin{equation}\label{E3-2}
\begin{aligned}
	&\left\langle \frac{\boldsymbol{u}_h^{k+\frac{1}{2}} - \boldsymbol{u}_h^{k-\frac{1}{2}}}{\varDelta t}, \boldsymbol{u}_h^{k+\frac{1}{2}} \right\rangle + \left\langle \boldsymbol{\omega}_h^{k-\frac{1}{2}} \times \boldsymbol{u}_h^{k+\frac{1}{2}}, \boldsymbol{u}_h^{k+\frac{1}{2}} \right\rangle + \left\langle \nabla \times \boldsymbol{\omega}_h^{k+\frac{1}{2}}, \boldsymbol{u}_h^{k+\frac{1}{2}} \right\rangle - \left\langle \phi_h^k, \nabla \cdot \boldsymbol{u}_h^{k+\frac{1}{2}} \right\rangle\\
	&\hspace{5.5cm}= - \left\langle \frac{p_h^{k-1} + p_h^{k}}{2} \nabla \mu_h^k, \boldsymbol{u}_h^{k+\frac{1}{2}} \right\rangle - \left\langle \frac{n_h^{k-1} + n_h^{k}}{2} \nabla \nu_h^k, \boldsymbol{u}_h^{k+\frac{1}{2}} \right\rangle.
\end{aligned}
\end{equation}
By utilizing the polarization identity, we can rewrite the first term in above equation as
\begin{equation*}
\left\langle \boldsymbol{u}_h^{k+\frac{1}{2}} - \boldsymbol{u}_h^{k-\frac{1}{2}}, \boldsymbol{u}_h^{k+\frac{1}{2}} \right\rangle = \frac{1}{2} \left( \left\| \boldsymbol{u}_h^{k+\frac{1}{2}} \right\|^2 - \left\| \boldsymbol{u}_h^{k-\frac{1}{2}} \right\|^2 + \left\| \boldsymbol{u}_h^{k+\frac{1}{2}} - \boldsymbol{u}_h^{k-\frac{1}{2}} \right\|^2 \right).
\end{equation*}
Next, the second term, namely the nonlinear convection term, is identically zero due to the fundamental vector identity $(\boldsymbol{a} \times \boldsymbol{b}) \cdot \boldsymbol{b} = 0$.
Additionally, by selecting the test function in \eqref{eq scheme 3b} to be $\boldsymbol{w}_h = \boldsymbol{\omega}_h^{k+\frac{1}{2}}$, the third term of \eqref{E3-2} can be rewritten as
\begin{equation*}
\left\langle \nabla \times \boldsymbol{\omega}_h^{k+\frac{1}{2}}, \boldsymbol{u}_h^{k+\frac{1}{2}} \right\rangle 
= \left\langle \boldsymbol{\omega}_h^{k+\frac{1}{2}}, \boldsymbol{\omega}_h^{k+\frac{1}{2}}\right\rangle 
= \left\| \boldsymbol{\omega}_h^{k+\frac{1}{2}} \right\|^2.
\end{equation*}
From \eqref{eq scheme 3c}, we know $\left\langle \phi_h^k, \nabla \cdot \boldsymbol{u}_h^{k+\frac{1}{2}} \right\rangle=0$ since \eqref{eq scheme 3c} holds for all $\alpha_h \in S(\Omega) $ and $ \phi_h^k\in S(\Omega)$.
Finally, \eqref{E3-2} becomes
\begin{equation*}
\begin{aligned}
	&\frac{1}{2\varDelta t} \left( \left\| \boldsymbol{u}_h^{k+\frac{1}{2}} \right\|^2 - \left\| \boldsymbol{u}_h^{k-\frac{1}{2}} \right\|^2 + \left\| \boldsymbol{u}_h^{k+\frac{1}{2}} - \boldsymbol{u}_h^{k-\frac{1}{2}} \right\|^2 \right) + \left\| \boldsymbol{\omega}_h^{k+\frac{1}{2}} \right\|^2 \\
	&\quad = - \left\langle \frac{p_h^{k-1} + p_h^{k}}{2} \nabla \mu_h^k, \boldsymbol{u}_h^{k+\frac{1}{2}} \right\rangle - \left\langle \frac{n_h^{k-1} + n_h^{k}}{2} \nabla \nu_h^k, \boldsymbol{u}_h^{k+\frac{1}{2}} \right\rangle.
\end{aligned}
\end{equation*}
Substituting $E_3$ and rearranging yields 
\begin{equation}\label{E3-final}
\begin{aligned}
	E_3 &= \frac{1}{2} \left( \left\| \boldsymbol{u}_h^{k+\frac{1}{2}} \right\|^2 - \left\| \boldsymbol{u}_h^{k-\frac{1}{2}} \right\|^2 \right)= - \varDelta t \left\| \boldsymbol{\omega}_h^{k+\frac{1}{2}}\right\|^2 - \varDelta t \left\langle \frac{p_h^{k-1} + p_h^{k}}{2} \nabla \mu_h^{k}, \boldsymbol{u}_h^{k+\frac{1}{2}}\right\rangle \\
	&\hspace{5.0cm} - \varDelta t \left\langle \frac{n_h^{k-1} + n_h^{k}}{2} \nabla \nu_h^{k}, \boldsymbol{u}_h^{k+\frac{1}{2}} \right\rangle - \frac{1}{2} \left\|\boldsymbol{u}_h^{k+\frac{1}{2}} - \boldsymbol{u}_h^{k-\frac{1}{2}} \right\|^2.
\end{aligned}
\end{equation}
From \eqref{eq E3} and \eqref{E3-final}, we obtain
\begin{equation}\label{Eq E}
\begin{split}
	\mathcal{E}_{h}^{k}-\mathcal{E}_{h}^{k-1} &= E_1 + E_2+E_3 \\
	&\leq -\frac{1}{2}\left\|\nabla \psi_h^k - \nabla \psi_h^{k-1}\right\|^2 
	- \varDelta t\left\|\boldsymbol{\omega}_h^{k+\frac{1}{2}}\right\|^2 
	- \frac{1}{2}\left\|\boldsymbol{u}_h^{k+\frac{1}{2}} - \boldsymbol{u}_h^{k-\frac{1}{2}}\right\|^2 \\
	&\quad - C'\left\|p_h^k - p_h^{k-1}\right\|^2 - C'\left\|n_h^k - n_h^{k-1}\right\|^2 \\
	&\quad + \varDelta t\left\langle \frac{p_h^{k-1} + p_h^{k}}{2}\nabla \mu_h^k, \boldsymbol{u}_h^{k-\frac{1}{2}} - \boldsymbol{u}_h^{k+\frac{1}{2}} \right\rangle 
	- \varDelta t\left\langle \frac{p_h^{k-1} + p_h^k}{2}\nabla \mu_h^{k-1}, \nabla \mu_h^k \right\rangle \\
	&\quad + \varDelta t\left\langle \frac{n_h^{k-1} + n_h^{k}}{2}\nabla \nu_h^k, \boldsymbol{u}_h^{k-\frac{1}{2}} - \boldsymbol{u}_h^{k+\frac{1}{2}} \right\rangle 
	- \varDelta t\left\langle \frac{n_h^{k-1} + n_h^k}{2}\nabla \nu_h^{k-1}, \nabla \nu_h^k \right\rangle.
\end{split}
\end{equation}
We now introduce notations
\[
\bar{p}_h^k = \frac{p_h^{k-1} + p_h^k}{2}, \quad \bar{n}_h^k = \frac{n_h^{k-1} + n_h^k}{2}.
\]
Then \eqref{Eq E} becomes
\begin{equation}\label{Eq E 1}
\begin{split}
	\mathcal{E}_{h}^{k}-\mathcal{E}_{h}^{k-1}&\leq -\frac{1}{2}\left\|\nabla \psi_h^k - \nabla \psi_h^{k-1}\right\|^2 
	- \varDelta t\left\|\boldsymbol{\omega}_h^{k+\frac{1}{2}}\right\|^2 
	- \frac{1}{2}\left\|\boldsymbol{u}_h^{k+\frac{1}{2}} - \boldsymbol{u}_h^{k-\frac{1}{2}}\right\|^2 \\
	&\quad - C'\left\|p_h^k - p_h^{k-1}\right\|^2 - C'\left\|n_h^k - n_h^{k-1}\right\|^2 \\
	&\quad + \varDelta t\left\langle \bar{p}_h^k \nabla \mu_h^k, \boldsymbol{u}_h^{k-\frac{1}{2}} - \boldsymbol{u}_h^{k+\frac{1}{2}} \right\rangle 
	- \varDelta t\left\langle \bar{p}_h^k \nabla \mu_h^{k-1}, \nabla \mu_h^k \right\rangle \\
	&\quad + \varDelta t\left\langle \bar{n}_h^k \nabla \nu_h^k, \boldsymbol{u}_h^{k-\frac{1}{2}} - \boldsymbol{u}_h^{k+\frac{1}{2}} \right\rangle 
	- \varDelta t\left\langle \bar{n}_h^k \nabla \nu_h^{k-1}, \nabla \nu_h^k \right\rangle.
\end{split}
\end{equation}
From Young's inequality, we know that
\[
\begin{aligned}
\varDelta t \left\langle \bar{p}_h^k \nabla \mu_h^k, \boldsymbol{u}_h^{k-\frac{1}{2}} - \boldsymbol{u}_h^{k+\frac{1}{2}} \right\rangle
&\leq \frac{1}{4}\left\| \boldsymbol{u}_h^{k+\frac{1}{2}} - \boldsymbol{u}_h^{k-\frac{1}{2}} \right\|^2 + \varDelta t^2\left\| \bar{p}_h^k \nabla \mu_h^k \right\|^2,\\
\varDelta t \left\langle \bar{n}_h^k \nabla \nu_h^k, \boldsymbol{u}_h^{k-\frac{1}{2}} - \boldsymbol{u}_h^{k+\frac{1}{2}} \right\rangle
&\leq \frac{1}{4}\left\| \boldsymbol{u}_h^{k+\frac{1}{2}} - \boldsymbol{u}_h^{k-\frac{1}{2}} \right\|^2 + \varDelta t^2\left\| \bar{n}_h^k \nabla \nu_h^k\right \|^2.
\end{aligned}
\]
Thus, \eqref{Eq E 1} can be rewrite as
\begin{equation}\label{Eq E 2}
\begin{split}
	\mathcal{E}_{h}^{k}-\mathcal{E}_{h}^{k-1} &\leq -\frac{1}{2}\left\|\nabla \psi_h^k - \nabla \psi_h^{k-1}\right\|^2 
	- \varDelta t\left\|\boldsymbol{\omega}_h^{k+\frac{1}{2}}\right\|^2  - C'\left\|p_h^k - p_h^{k-1}\right\|^2 - C'\left\|n_h^k - n_h^{k-1}\right\|^2 \\
	&\quad  + \varDelta t^2\left\| \bar{p}_h^k \nabla \mu_h^k \right\|^2
	- \varDelta t\left\langle \bar{p}_h^k\nabla \mu_h^{k-1}, \nabla \mu_h^k \right\rangle + \varDelta t^2\left\| \bar{n}_h^k \nabla \nu_h^k \right\|^2 
	- \varDelta t\left\langle \bar{n}_h^k\nabla \nu_h^{k-1}, \nabla \nu_h^k \right\rangle.
\end{split}
\end{equation}
We now estimate the term
$
- \varDelta t \left\langle \bar{p}_h^k \nabla \mu_h^{k-1}, \nabla \mu_h^k \right\rangle 
$
and use the algebraic identity for arbitrary vectors $a$ and $b$,
$
- a b = \frac{1}{2} |a - b|^2 - \frac{1}{2} |a|^2 - \frac{1}{2} |b|^2.
$ We can obtain
\begin{equation*}
\begin{aligned}
	- \varDelta t \left\langle \bar{p}_h^k \nabla \mu_h^{k-1}, \nabla \mu_h^k \right\rangle & = \frac{\varDelta t}{2} \int_{\Omega}\bar{p}_h^k  \left( \left| \nabla \mu_h^k - \nabla \mu_h^{k-1} \right|^2 - \left| \nabla \mu_h^k \right|^2 - \left| \nabla \mu_h^{k-1} \right|^2 \right) \mathrm{d}x\\ &\leq \frac{\varDelta t}{2} M \left\| \nabla \mu_h^k - \nabla \mu_h^{k-1} \right\|^2 - \frac{\varDelta t}{2} \left\| \sqrt{\bar{p}_h^k } \nabla \mu_h^k \right\|^2 - \frac{\varDelta t}{2} \left\| \sqrt{\bar{p}_h^k } \nabla \mu_h^{k-1} \right\|^2.  	
\end{aligned}
\end{equation*}
Next, we estimate 
\begin{equation*}
\begin{aligned}
	\varDelta t^2 \left\Vert \bar{p}_h^k \nabla \mu_h^k \right\Vert^2- \frac{\varDelta t}{2} \left\| \sqrt{\bar{p}_h^k } \nabla \mu_h^k \right\|^2 
	&= \int_{\Omega} \left( \varDelta t^2 \left(\bar{p}_h^k\right)^2 - \frac{\varDelta t}{2} \bar{p}_h^k \right) \left|\nabla \mu_h^k\right|^2 \, \mathrm{d}x\\& \leq
	\varDelta t \left( \varDelta t M^2 - \frac{m}{2} \right) \left\|\nabla \mu_h^k\right\|^2 .
\end{aligned}
\end{equation*}
Combining these two relations, we end up with
\begin{equation}\label{E3 a1}
\begin{aligned}
	& \varDelta t^2 \left\Vert \bar{p}_h^k \nabla \mu_h^k \right\Vert^2 - \varDelta t \left\langle \bar{p}_h^k \nabla \mu_h^{k-1}, \nabla u_h^k \right\rangle \\
	&\quad\leq \varDelta t \left( \varDelta t M^2 - \frac{m}{2} \right) \left\|\nabla \mu_h^k\right\|^2  + \frac{\varDelta t}{2} M \left \Vert   \nabla \mu_h^k - \nabla \mu_h^{k-1}  \right\Vert^2  - \frac{\varDelta t}{2} M \left\Vert  \nabla \mu_h^{k-1} \right\Vert^2.
\end{aligned}
\end{equation}
Repeating the derivation for the term involving $\nu_h$, we get
\begin{equation}\label{E3 a2}
\begin{aligned}
	& \varDelta t^2 \left\Vert\bar{n}_h^k \nabla \nu_h^k \right\Vert^2 - \varDelta t \left\langle \bar{n}_h^k \nabla \nu_h^{k-1}, \nabla \nu_h^k \right\rangle \\
	&\quad\leq  \varDelta t \left( \varDelta t M^2 - \frac{m}{2} \right) \left\|\nabla \nu_h^k\right\|^2 + \frac{\varDelta t}{2} M \left\Vert \nabla \nu_h^k - \nabla \nu_h^{k-1} \right\Vert^2 - \frac{\varDelta t}{2} M \left\Vert \nabla \nu_h^{k-1} \right\Vert^2.
\end{aligned}
\end{equation}
Inserting \eqref{E3 a1} and \eqref{E3 a2} into \eqref{Eq E 2} leads to
\begin{equation}\label{Eq E 4}
\begin{split}
	\mathcal{E}_{h}^{k}-\mathcal{E}_{h}^{k-1} &\leq -\frac{1}{2}\left\|\nabla \psi_h^k - \nabla \psi_h^{k-1}\right\|^2 
	- \varDelta t\left\|\boldsymbol{\omega}_h^{k+\frac{1}{2}}\right\|^2  - C'\left\|p_h^k - p_h^{k-1}\right\|^2 - C'\left\|n_h^k - n_h^{k-1}\right\|^2 \\
	&\quad  +\varDelta t \left( \varDelta t M^2 - \frac{m}{2} \right) \left\|\nabla \mu_h^k\right\|^2  + \frac{\varDelta t}{2} M \left \Vert   \nabla \mu_h^k - \nabla \mu_h^{k-1}  \right\Vert^2  - \frac{\varDelta t}{2} M \left\Vert  \nabla \mu_h^{k-1} \right\Vert^2 \\
	&\quad+ \varDelta t \left( \varDelta t M^2 - \frac{m}{2} \right) \left\|\nabla \nu_h^k\right\|^2  + \frac{\varDelta t}{2} M \left\Vert \nabla \nu_h^k - \nabla \nu_h^{k-1} \right\Vert^2 - \frac{\varDelta t}{2} M \left\Vert  \nabla \nu_h^{k-1} \right\Vert^2.
\end{split}
\end{equation}
Under the constraint $\varDelta t <  \frac{m}{2M^2} $, the terms $\varDelta t \left( \varDelta t M^2 - \frac{m}{2} \right) \left\|\nabla \mu_h^k\right\|^2 $ are strictly negative, and the terms $- \dfrac{\varDelta t}{2} M \left\Vert  \nabla \mu_h^{k-1} \right\Vert^2$ and $- \dfrac{\varDelta t}{2} M \left\Vert  \nabla \nu_h^{k-1} \right\Vert^2$ are non-positive. Therefore, we only need to further prove that
\begin{equation}\label{F}
\begin{split}
	&F: = - \frac{1}{2} \left\Vert \nabla \psi_h^k -	 \nabla \psi_h^{k-1} \right\Vert^2 + \frac{\varDelta t}{2} M \left\Vert \nabla \mu_h^k - \nabla \mu_h^{k-1} \right\Vert^2 \\&\hspace{1cm}+ \frac{\varDelta t}{2} M \left\Vert \nabla \nu_h^k - \nabla \nu_h^{k-1} \right\Vert^2 - C' \left\Vert p_h^k - p_h^{k-1} \right\Vert^2 - C' \left\Vert n_h^k - n_h^{k-1} \right\Vert^2
\end{split}
\end{equation}
is negative.
It follows from Lemma 2 and Equation (A.1) in the Appendix that
\begin{equation*}
\frac{\varDelta t}{2} M \left\|\nabla \mu_h^k - \nabla \mu_h^{k-1}\right\|^2\leq C \varDelta t h^{-2}\left\| p_h^k - p_h^{k-1} \right\|^2+ C \varDelta t h^{-2} \left\|\nabla \psi_h^k - \nabla \psi_h^{k-1}\right\|^2,
\end{equation*}
where $h$ is a constant related to the mesh size. Similarly, we can obtain
\begin{equation*}
\frac{\varDelta t}{2} M\left\|  \nabla \nu_h^k - \nabla \nu_h^{k-1}  \right\|^2 \leq C \varDelta t h^{-2}\left\| n_h^k - n_h^{k-1} \right\|^2 +  C \varDelta t h^{-2} \left\|\nabla \psi_h^k - \nabla \psi_h^{k-1}\right\|^2.
\end{equation*}
Substituting these relations into \eqref{F}, we have
\begin{equation*}
\begin{aligned}
	F &\leq  \left( -\frac{1}{2} +  C \varDelta t h^{-2} \right) \left\|\nabla \psi_h^k - \nabla \psi_h^{k-1}\right\|^2 + \left (C \varDelta t h^{-2}- C' \right)\left\| p_h^k - p_h^{k-1} \right\|^2  \\
	&\qquad\qquad \qquad+\left (C \varDelta t h^{-2}- C' \right)\left\| n_h^k - n_h^{k-1} \right\|^2 .
\end{aligned}
\end{equation*}
This eventually gives
\begin{equation}
\begin{split}
	\mathcal{E}_{h}^{k}-\mathcal{E}_{h}^{k-1} 	&\leq \left( -\frac{1}{2} +  C \varDelta t h^{-2}\right) \left\|\nabla \psi_h^k - \nabla \psi_h^{k-1}\right\|^2
	- \varDelta t\left\|\boldsymbol{\omega}_h^{k+\frac{1}{2}}\right\|^2  \\&\quad- \frac{\varDelta t}{2} M \left\Vert  \nabla \mu_h^{k-1} \right\Vert^2- \frac{\varDelta t}{2} M \left\Vert  \nabla \nu_h^{k-1} \right\Vert^2  \\
	&\quad+\varDelta t \left( \varDelta t M^2 - \frac{m}{2} \right) \left\|\nabla \mu_h^k\right\|^2  + \varDelta t \left( \varDelta t M^2 - \frac{m}{2} \right) \left\|\nabla \nu_h^k\right\|^2   \\
	&\quad+ \left(C  \varDelta t h^{-2}-C'\right) \left\|p_h^k - p_h^{k-1}\right\|^2 + \left(C \varDelta t h^{-2}-C'\right)\left\|n_h^k - n_h^{k-1}\right\|^2,
\end{split}
\end{equation}
which implies the stability conditions 
\[
-\frac{1}{2} +  C \varDelta t h^{-2} < 0,\quad\varDelta t M^2 - \frac{m}{2}<0
, \quad \text{and} \quad C \varDelta t h^{-2} - C' < 0.
\]
Finally, we conclude that the stability of the numerical scheme,
$$\mathcal{E}_{h}^{k}-\mathcal{E}_{h}^{k-1} < 0,$$
is guaranteed provided that the time step satisfies
$
\varDelta t \leq  \min\left\{ \frac{m}{2M^2}, {C h^2} \right\},
$	where the constant $C$ depends on $M$, $m$, $C_{\mathrm{inv}}$, and $C_P$ from Lemma 2 in the Appendix, but remains independent of the mesh size $h$.
The proof of is complete.

\section{Numerical tests}
\label{SEC: numeric}

In this section, we numerically test the accuracy, the mass conservation and the conditional energy stability of the present method using manufactured solutions and a benchmark problem. The finite-dimensional spaces used here are the mimetic spectral element spaces \cite{zhang2022mimetic}. Any other set of finite-dimensional spaces that satisfies the discrete de Rham complex \eqref{Eq: de} and the additional regularity works as well.

\subsection{Accuracy tests}
We construct manufactured solutions  of $p$, $\psi$, $\boldsymbol{u}$ and $P$,
\begin{subequations}\label{Eq: manufactured solutions}
\begin{align*}
	p &= \left[\cos(x)\sin(y)+3\right]e^{t},\\
	\psi & = \sin(x)\cos(y) e^{t},\\
	\boldsymbol{u} &= \begin{bmatrix}
		\sin(x)\cos(y)e^t&
		-\cos(x)\sin(y)e^t
	\end{bmatrix}^{\mathsf{T}},\\
	P & = \sin(x)\sin(y) e^t + \frac{1}{2}\boldsymbol{u}\cdot\boldsymbol{u}.
\end{align*}
\end{subequations}
Manufactured solutions of  $\phi$, $n$ and  $\omega$  then can be computed through \eqref{Eq: phi}, \eqref{Eq: PNPNSc} and \eqref{Eq: PNPNSe}, respectively. Note that extra source terms need to be added to evolution equations \eqref{Eq: PNPNSa}, \eqref{Eq: PNPNSb}, and \eqref{Eq: PNPNSd} to balance them under these manufactured solutions. And the weak formulation and its discretization need to be modified accordingly. This is common for PDEs that are hard to find a set of analytical solutions. 

The computational domain is selected to be $\Omega=(0, 2\pi)^2$ and is fully periodic, i.e., $\partial\Omega=\emptyset$. Uniform meshes of $K\times K$ square elements are generated in $\Omega$, and we use $N$ to denote the degree of the finite-dimensional spaces. The time step size is denoted by $\varDelta t$. Provided with initial conditions at $t^0=0$ and source terms according to the manufactured solutions, we compute the PNPNS problem using the present method. When $t^k=1$, errors
\[
\left\|p^{k}_h\right\|_{L^2\text{-error}},
\left\|n^{k}_h\right\|_{L^2\text{-error}},
\left\|\psi^{k}_h\right\|_{H^1\text{-error}},
\left\|\boldsymbol{u}^{k+\frac{1}{2}}_h\right\|_{H(\mathrm{div})\text{-error}},
\left\|\omega^{k+\frac{1}{2}}_h\right\|_{H(\mathrm{curl})\text{-error}}\ \text{and}\ 
\left\|\phi^{k}_h\right\|_{L^2\text{-error}}
\]
are measured, where, for example, 
\[\left\|p^{k}_h\right\|_{L^2\text{-error}} := \left\|p^{k}_h - p(t^k, \boldsymbol{x})\right\|_{L^2}\ \text{and}\  \left\|\psi^{k}_h\right\|_{H^1\text{-error}} := \sqrt{\left\|\psi^{k}_h\right\|^2_{L^2\text{-error}} +  \left\|\nabla \psi^{k}_h\right\|^2_{L^2\text{-error}}}.\]

\subsubsection{Temporal convergence rates} \label{SUB: tcr}
\begin{figure}[h!]
\centering
\begin{minipage}[c]{1\textwidth}
	\centering{
		\subfloat{ \label{TC a}
			\begin{minipage}[b]{0.33\textwidth}
				\centering
				\includegraphics[width=1\linewidth]{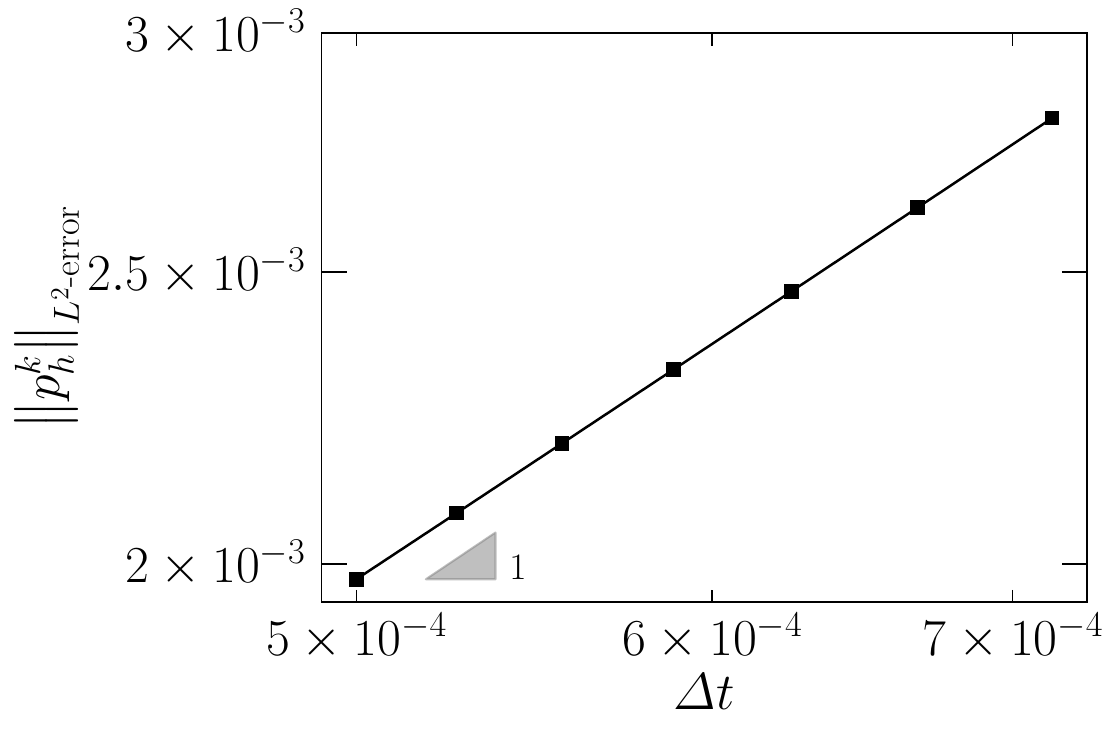}
			\end{minipage}
		}
		\subfloat{ \label{TC b}
			\begin{minipage}[b]{0.33\textwidth}
				\centering
				\includegraphics[width=1\linewidth]{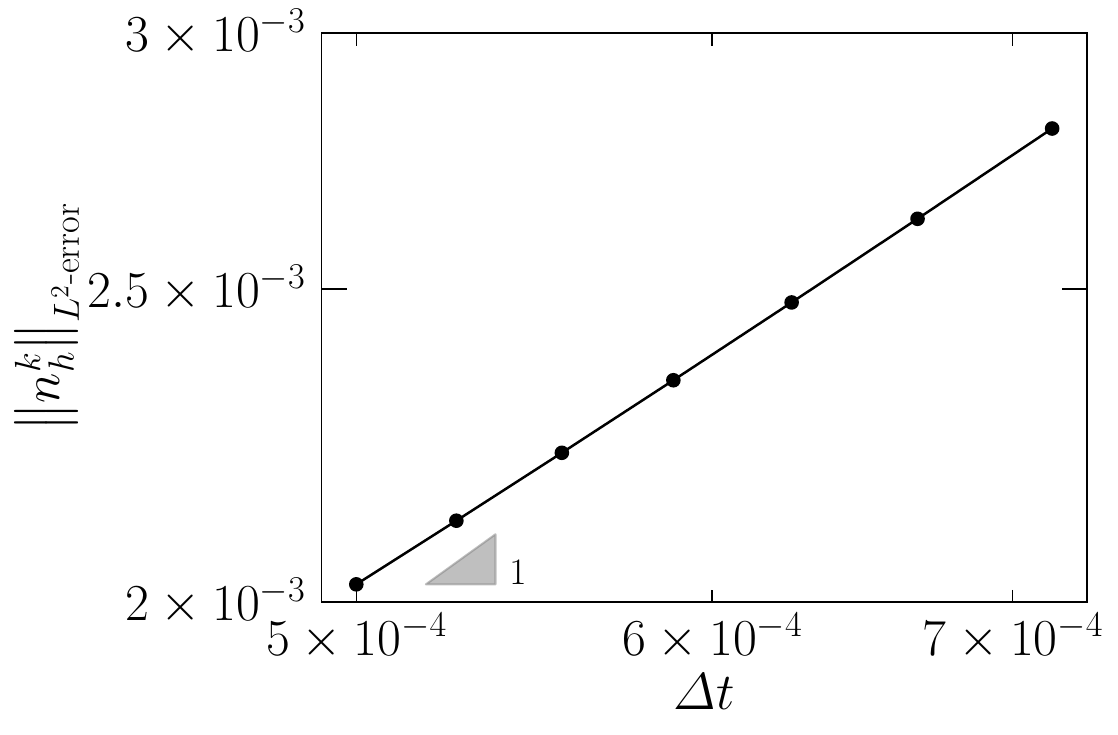}
			\end{minipage}
		}
		\subfloat{ \label{TC c}
			\begin{minipage}[b]{0.33\textwidth}
				\centering
				\includegraphics[width=1\linewidth]{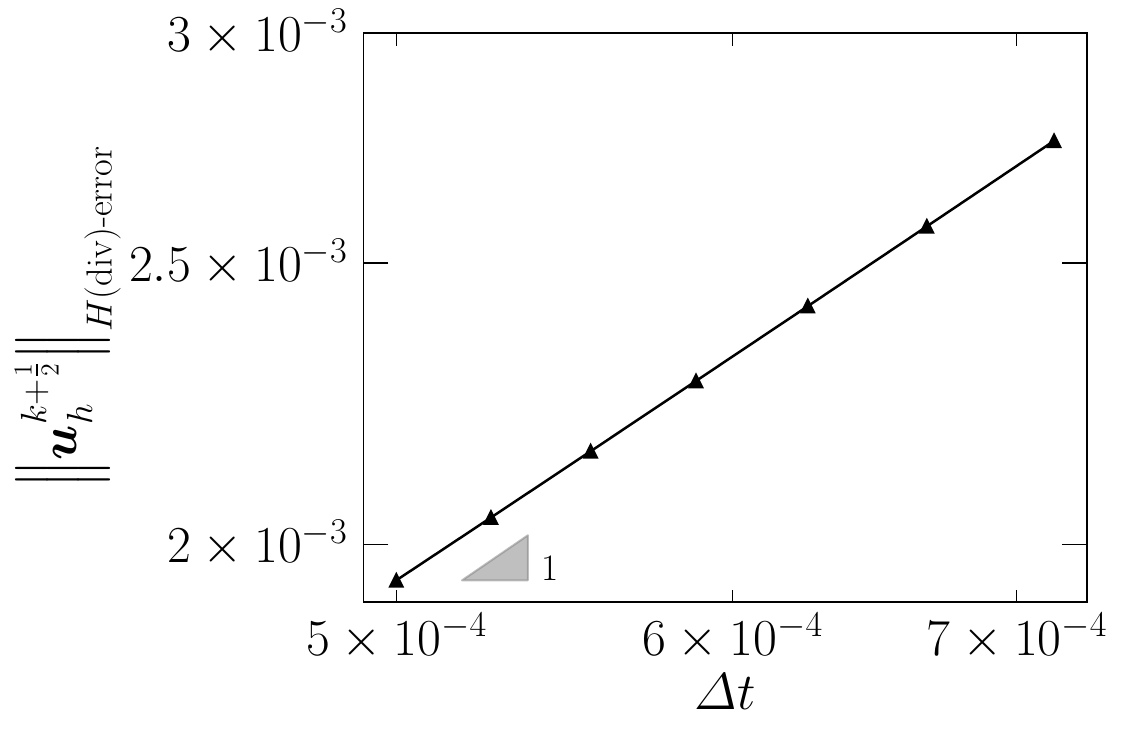}
			\end{minipage}
		}
	}
\end{minipage}%
\caption{Some results of the temporal accuracy tests using $K=30$, $N=3$ and $\varDelta t \in \left\lbrace\frac{1}{1400}, \frac{1}{1500}, \cdots, \frac{1}{2000}\right\rbrace$.} 
\label{fig:TC}
\end{figure}

To examine the temporal accuracy, we run tests using $K=30$, $N=3$ and $\varDelta t \in \left\lbrace \frac{1}{1400}, \frac{1}{1500}, \cdots, \frac{1}{2000}\right\rbrace$ and present some results in Fig.~\ref{fig:TC}. These results show that the proposed method has a first-order temporal accuracy which is inline with fact that an implicit-explicit temporal treatment is employed.

\subsubsection{Spatial convergence rates}
\begin{figure}[h!]
\centering
\begin{minipage}[c]{1\textwidth}
	\centering{
		\subfloat{ \label{SC a}
			\begin{minipage}[b]{0.33\textwidth}
				\centering
				\includegraphics[width=1\linewidth]{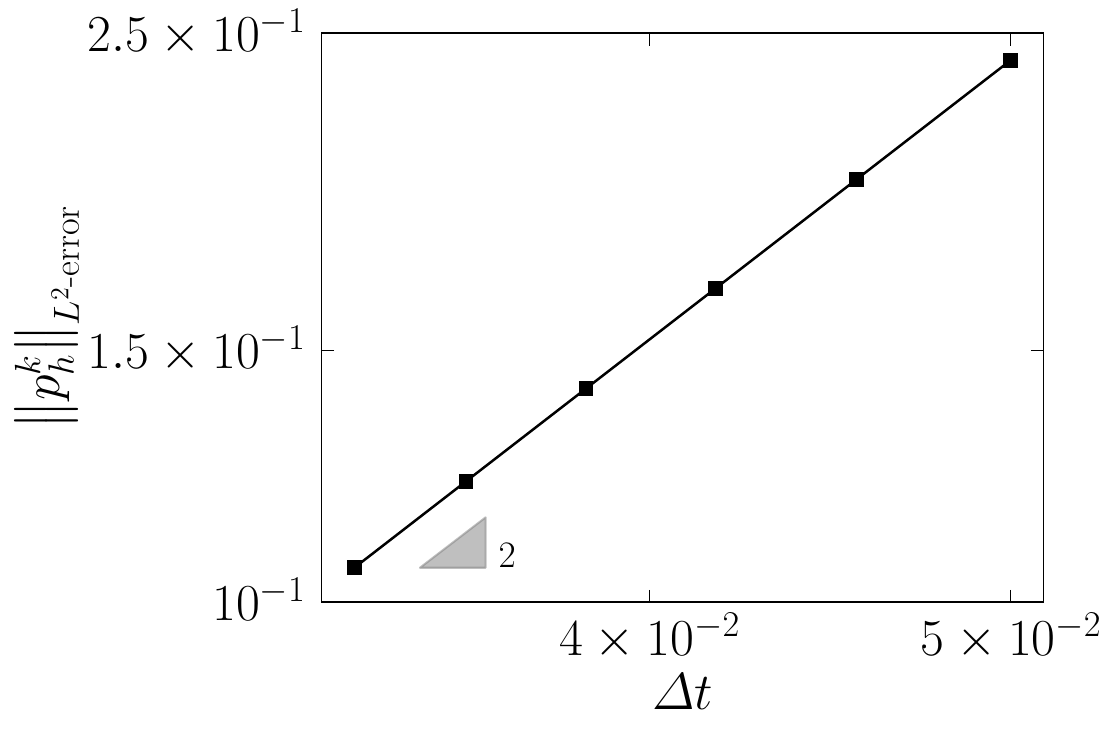}
			\end{minipage}
		}
		\subfloat{ \label{SC b}
			\begin{minipage}[b]{0.33\textwidth}
				\centering
				\includegraphics[width=1\linewidth]{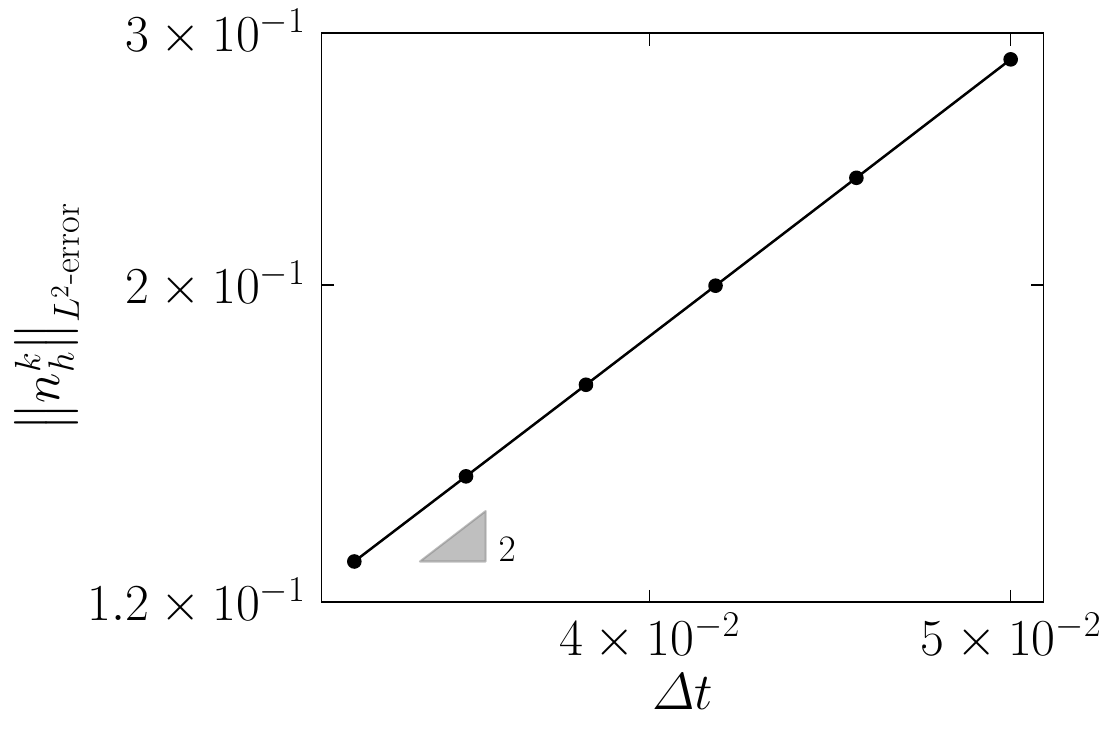}
			\end{minipage}
		}
		\subfloat{ \label{SC c}
			\begin{minipage}[b]{0.33\textwidth}
				\centering
				\includegraphics[width=1\linewidth]{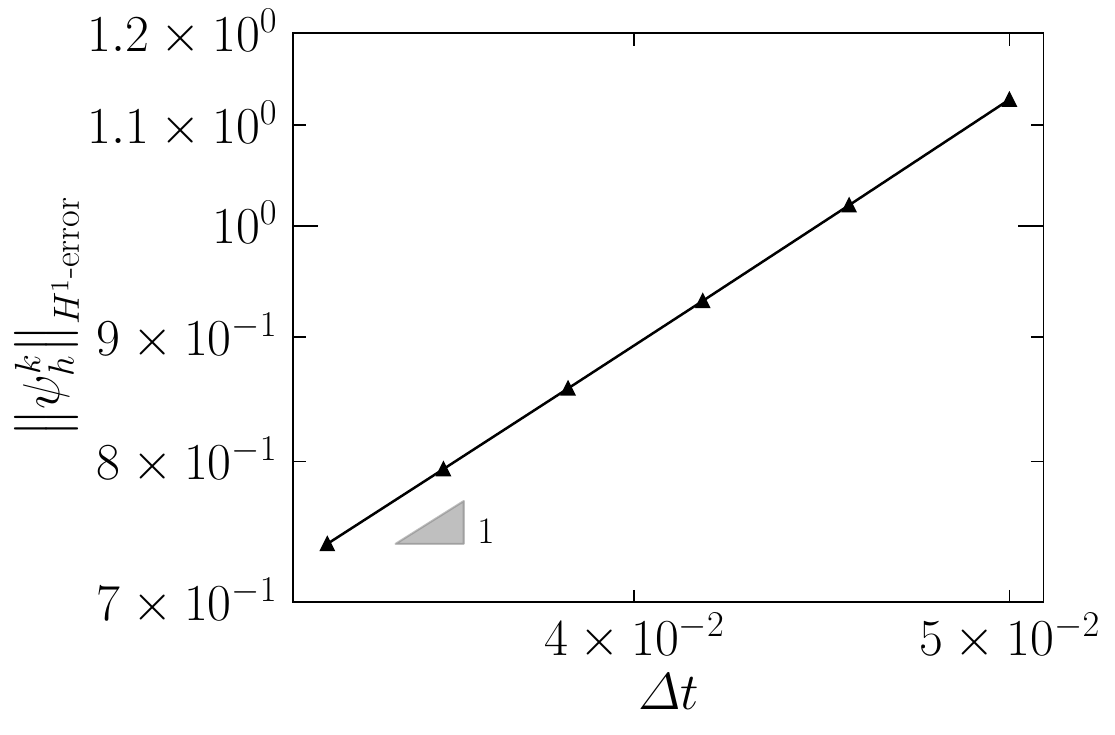}
			\end{minipage}
		}\\
		\subfloat{ \label{SC d}
			\begin{minipage}[b]{0.33\textwidth}
				\centering
				\includegraphics[width=1\linewidth]{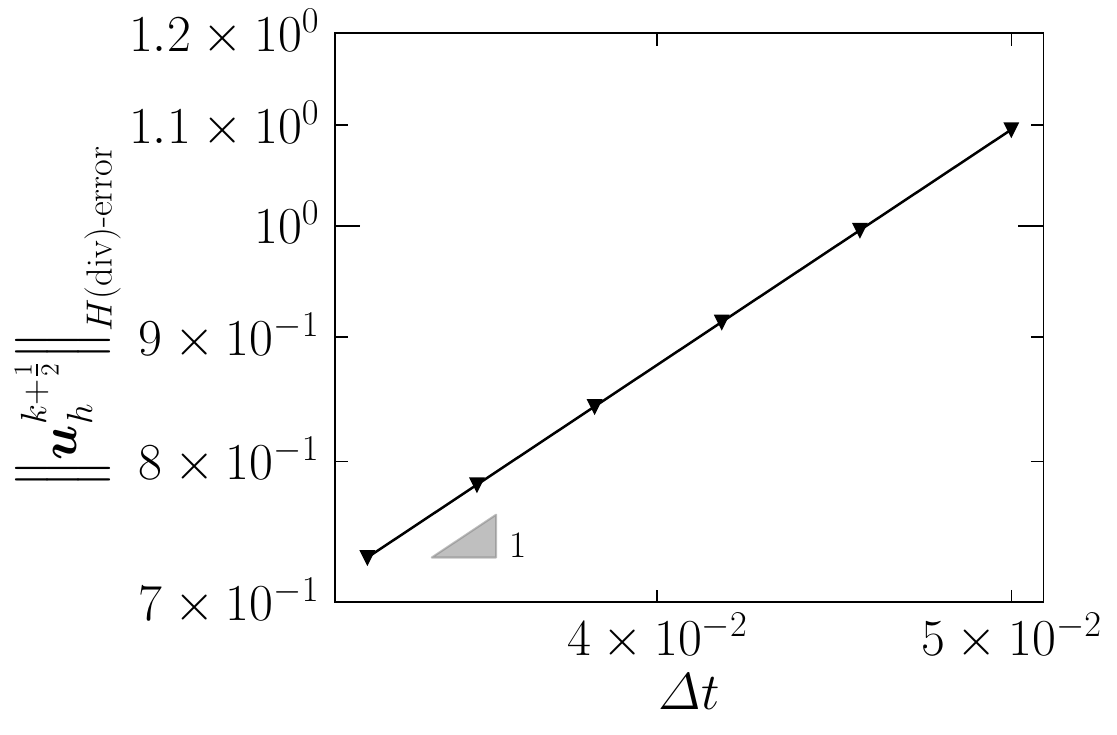}
			\end{minipage}
		}
		\subfloat{ \label{SC e}
			\begin{minipage}[b]{0.33\textwidth}
				\centering
				\includegraphics[width=1\linewidth]{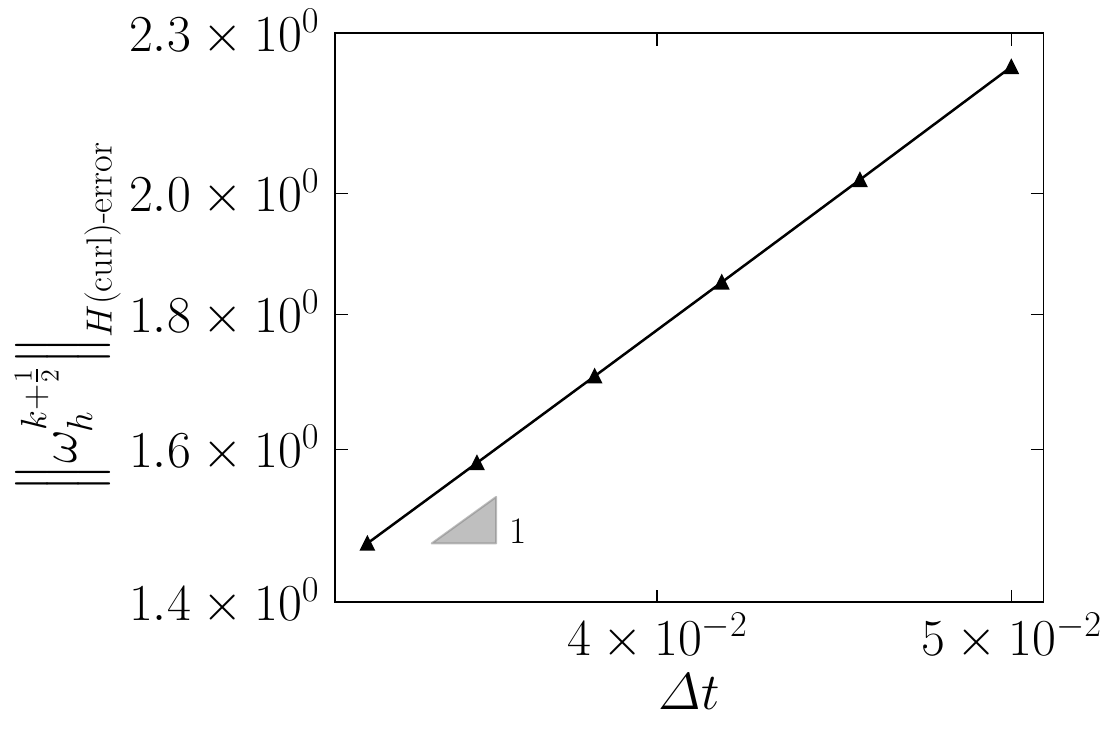}
			\end{minipage}
		}
		\subfloat{ \label{SC f}
			\begin{minipage}[b]{0.33\textwidth}
				\centering
				\includegraphics[width=1\linewidth]{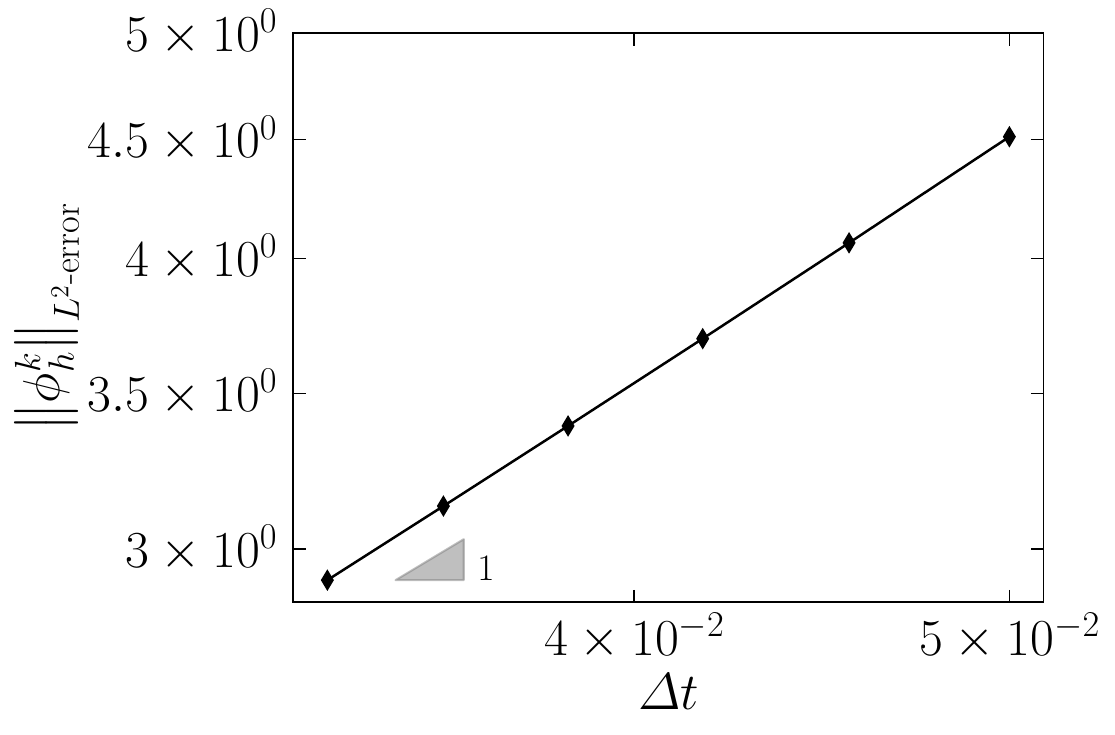}
			\end{minipage}
		}
	}
\end{minipage}%
\caption{Some results of the spatial accuracy tests under an $h$-refinement for $N=1$, $\varDelta t=\frac{1}{2000}$ and $K \in \left\lbrace20, 22,\cdots , 30\right\rbrace$.} 
\label{fig:SC}
\end{figure}

In Fig.~\ref{fig:SC}, results of an $h$-refinement for $N=1$, $K\in\left\lbrace20,22,\cdots,30\right\rbrace$ are presented where, to eliminate the pollution of temporal error, we have used $\varDelta t = \frac{1}{2000}$ (cf. Section~\ref{SUB: tcr}). These results imply that the proposed method has optimal spatial convergence rates for all unknowns. 

\subsection{A test for the mass conservation and the  energy stability}
In this test, the domain is selected to be $\Omega=(-1, 1)^2$. A set of manufactured initial conditions is given as
\begin{subequations}\label{Eq: IC}
\begin{align*}
	p^{0} &= 1.1 + \cos(\pi x)\cos(\pi y),\\
	n^{0} &= 1.1 - \cos(\pi x)\cos(\pi y),\\
	\boldsymbol{u}^0 & = \begin{bmatrix}
		\pi  \sin(\pi x) ^2 \sin(2\pi y) & - \pi \sin(2\pi x)  \sin(\pi y) ^2
	\end{bmatrix}^{\mathsf{T}},
\end{align*}
\end{subequations}
and $\omega^{0} =\nabla\times\boldsymbol{u}^0$. Suppose $\psi^0$ satisfies an Neumann boundary condition $\nabla\psi^0\cdot\boldsymbol{n} = 0$. We can easily find the analytical expression of $\psi^0$,
\[
\psi^{0} = \dfrac{1}{\epsilon\pi^2}\cos(\pi x)\cos(\pi y).
\]
by solving the Poisson equation \eqref{Eq: PNPNSc}. One can provide that, in $\Omega=(-1, 1)^2$, this set of initial conditions also satisfies boundary conditions \eqref{Eq: boundary conditions}. 
With these initial conditions and boundary conditions \eqref{Eq: boundary conditions}, the PNPNS problem is solved until it reaches $t^k=1$ using the proposed method for $N=2$, $\varDelta t = \frac{1}{2000}$ on a uniform mesh of $20\times 20$ square elements in $\Omega$. 

\begin{figure}[h!]
\centering
\begin{minipage}[c]{1\textwidth}
	\centering{
		\subfloat{ 
			\begin{minipage}[b]{0.33\textwidth}
				\centering
				\includegraphics[width=1\linewidth]{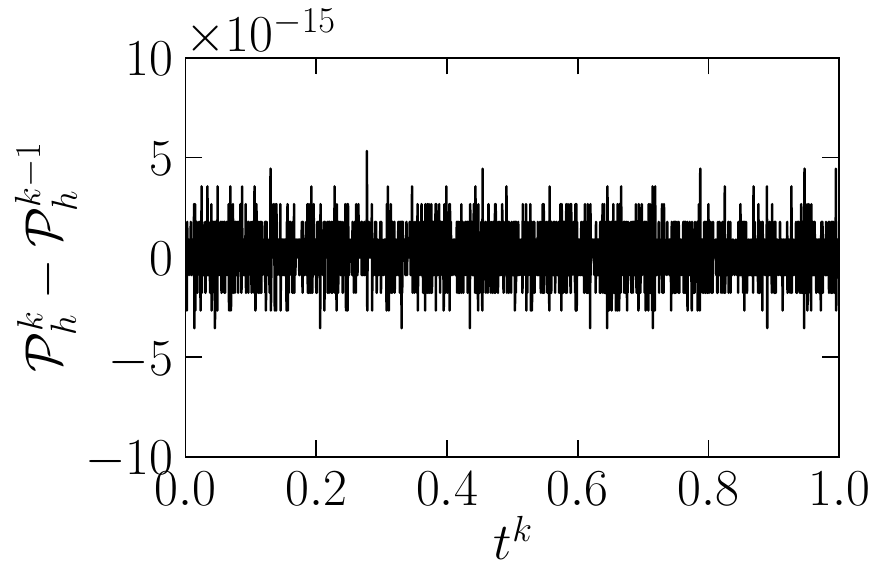}
			\end{minipage}
		}
		\subfloat{
			\begin{minipage}[b]{0.33\textwidth}
				\centering
				\includegraphics[width=1\linewidth]{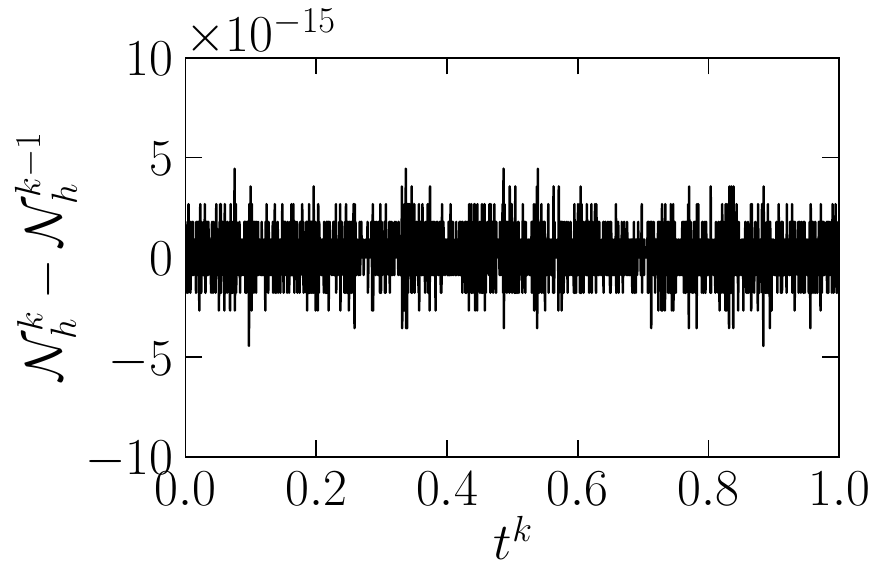}
			\end{minipage}
		}
		\subfloat{
			\begin{minipage}[b]{0.31\textwidth}
				\centering
				\includegraphics[width=1\linewidth]{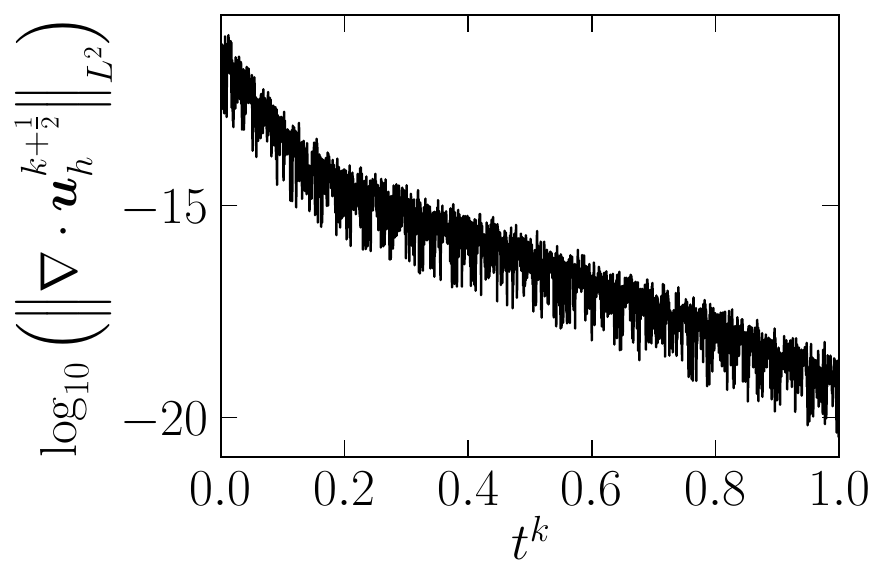}
			\end{minipage}
		}
	}
\end{minipage}%
\caption{Mass conservation results of the present method.} 
\label{fig: mass conservation}
\end{figure}

\begin{figure}[h!]
\centering
\begin{minipage}[c]{1\textwidth}
	\centering{
		\subfloat{ 
			\begin{minipage}[b]{0.35\textwidth}
				\centering
				\includegraphics[width=0.9\linewidth]{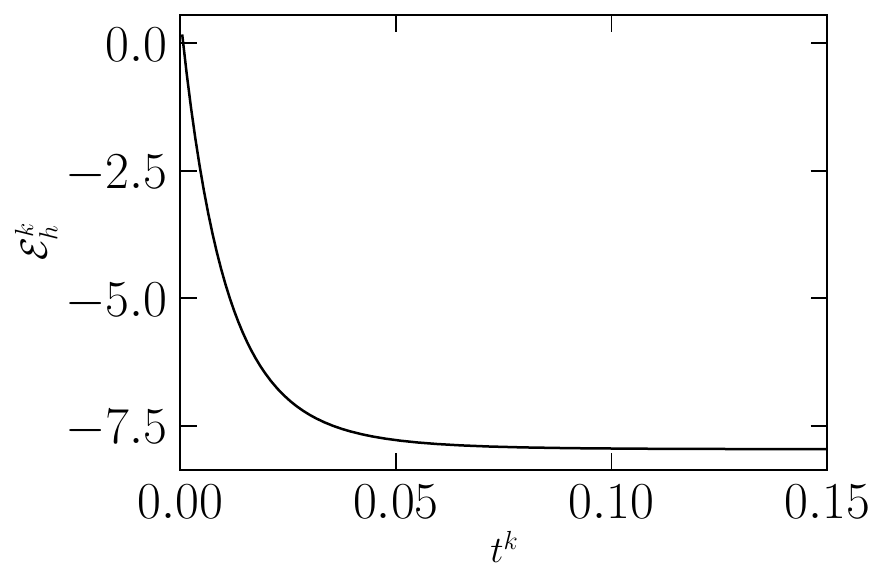}
			\end{minipage}
		}
		\subfloat{
			\begin{minipage}[b]{0.36\textwidth}
				\centering
				\includegraphics[width=0.9\linewidth]{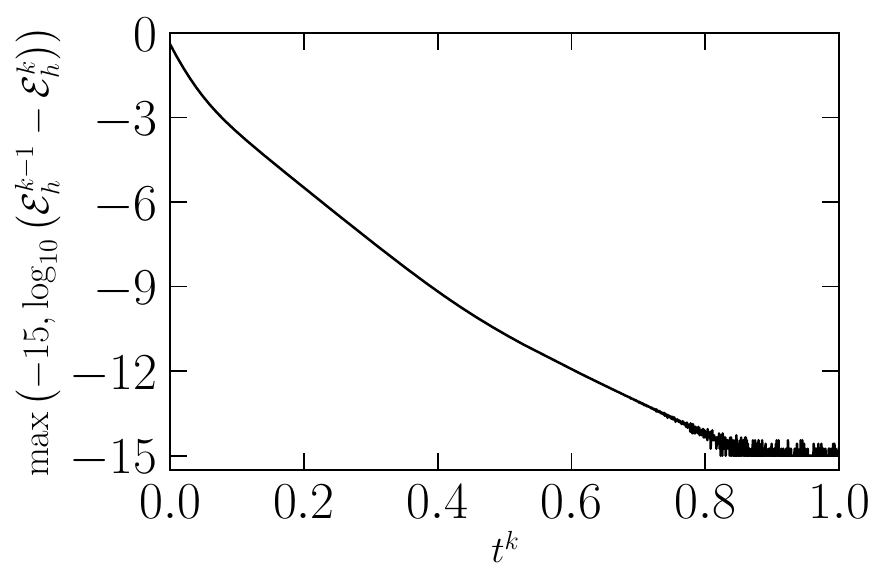}
			\end{minipage}
		}
	}
\end{minipage}%
\caption{Energy conservation results of the present method. Note that we assume the flow reaches an energy-conserving state if $ -10^{-15} < \mathcal{E}_{h}^{k-1}-\mathcal{E}_{h}^{k} < 0 $. And in the right diagram only the results for $ \log_{10}\left( \mathcal{E}_{h}^{k-1}-\mathcal{E}_{h}^{k}\right) \geq -15 $ are plotted.} 
\label{fig: EC}
\end{figure}

In Fig.~\ref{fig: mass conservation}, results showing that the mass conservation is exactly satisfied up to a machine precision are present. And in Fig.~\ref{fig: EC}, the energy stability is demonstrated. When dissipations, for example, the viscous dissipation, of the flow are strong, say when $t^k < 0.05$, the results are not convincing because, even if a significant amount of artificial energy is generated by the numerical scheme, the strong dissipations will overwhelm the artificial energy. But, when $t^k > 0.5$, the dissipations become extremely weak, and the present method still enjoys an  energy stability. 

\subsection{A benchmark test of initial discontinuous concentrations}
\begin{figure}[!htb]
\centering
\begin{minipage}[c]{0.9\textwidth}
	\centering{
		\subfloat{
			\begin{minipage}[b]{0.235\textwidth}
				\centering
				\includegraphics[width=1\linewidth]{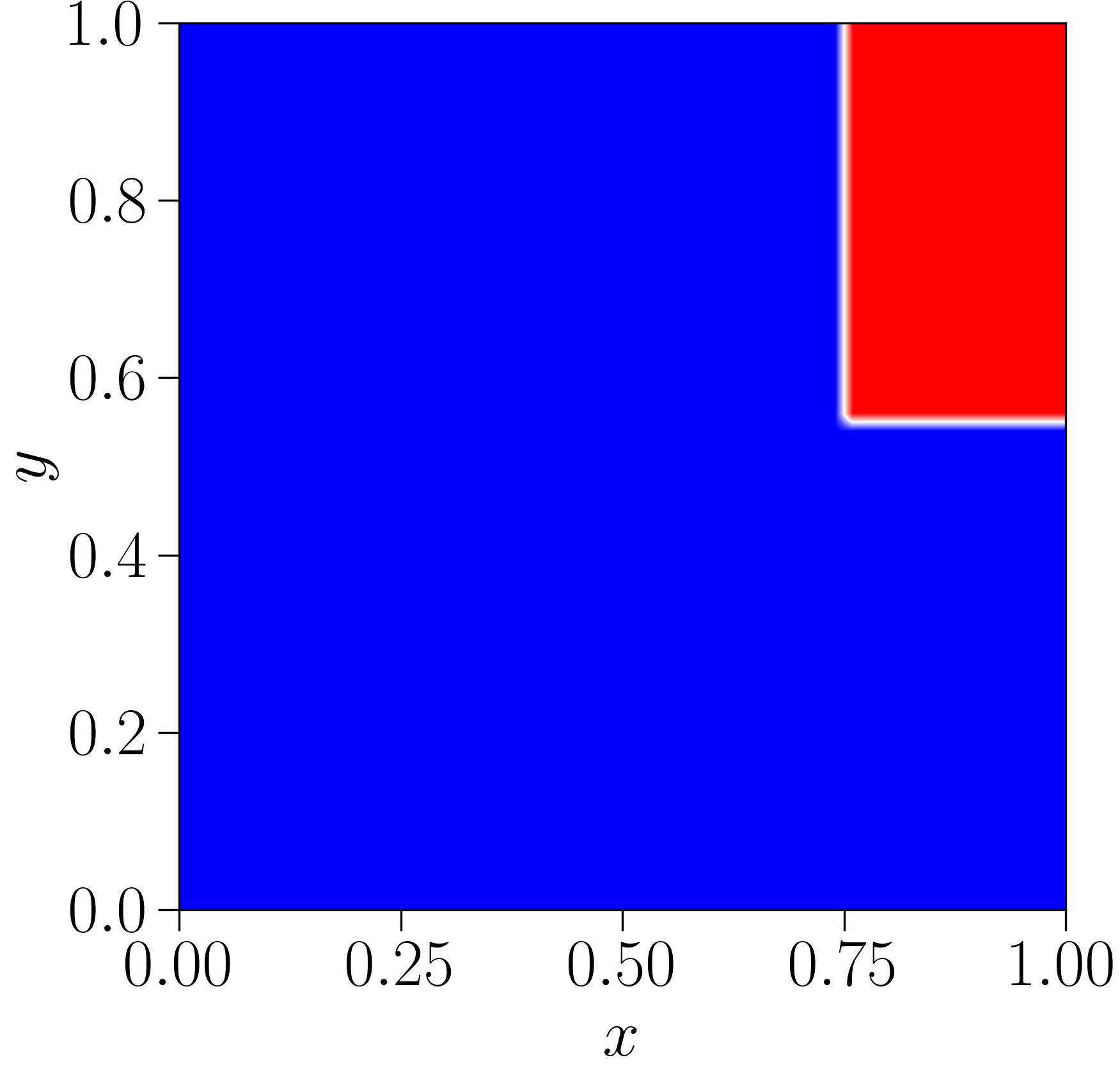}
			\end{minipage}
		}
		\subfloat{
			\begin{minipage}[b]{0.235\textwidth}
				\centering
				\includegraphics[width=1\linewidth]{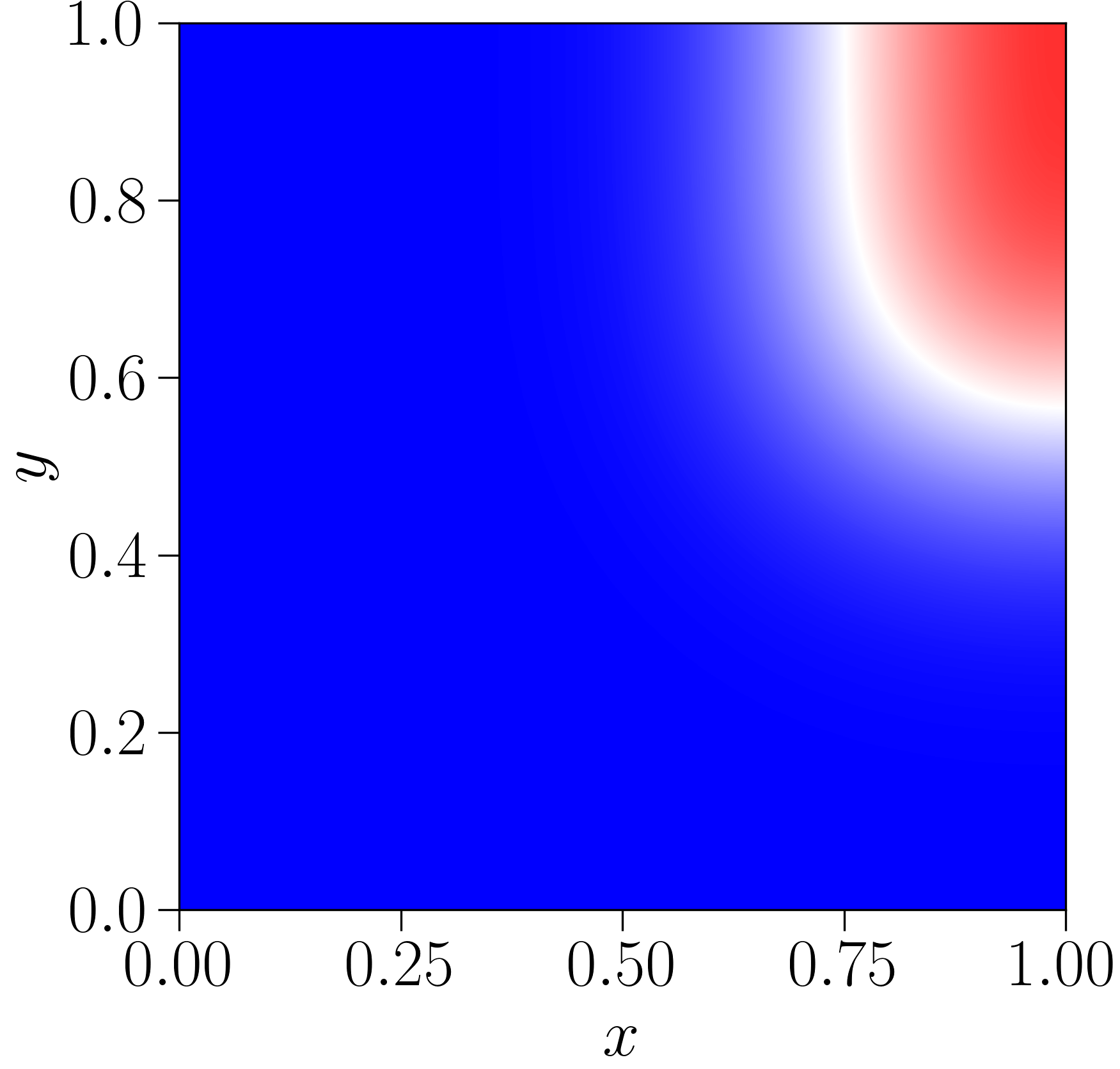}
			\end{minipage}
		}
		\subfloat{
			\begin{minipage}[b]{0.235\textwidth}
				\centering
				\includegraphics[width=1\linewidth]{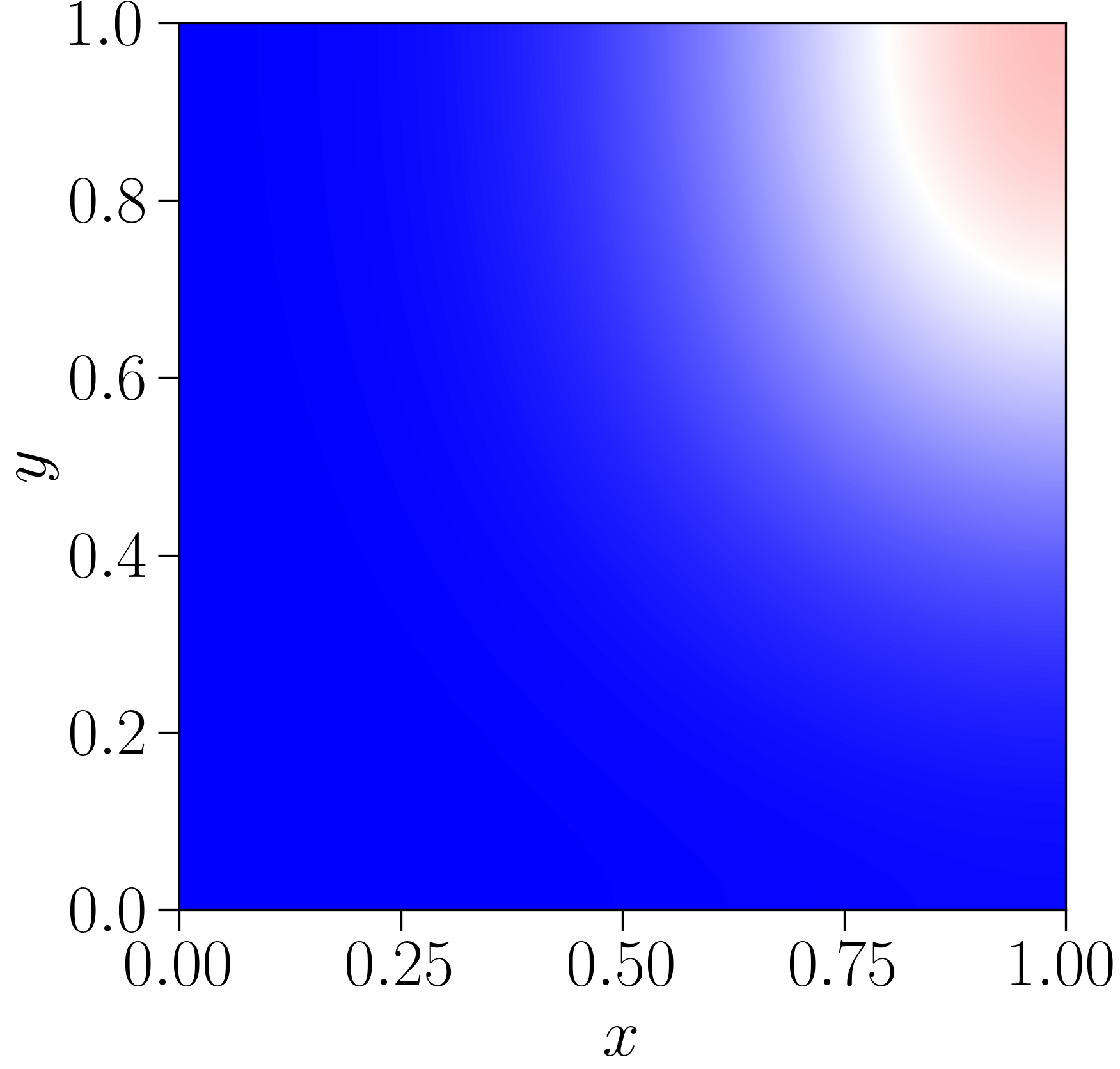}
			\end{minipage}
		}
		\subfloat{
			\begin{minipage}[b]{0.235\textwidth}
				\centering
				\includegraphics[width=1\linewidth]{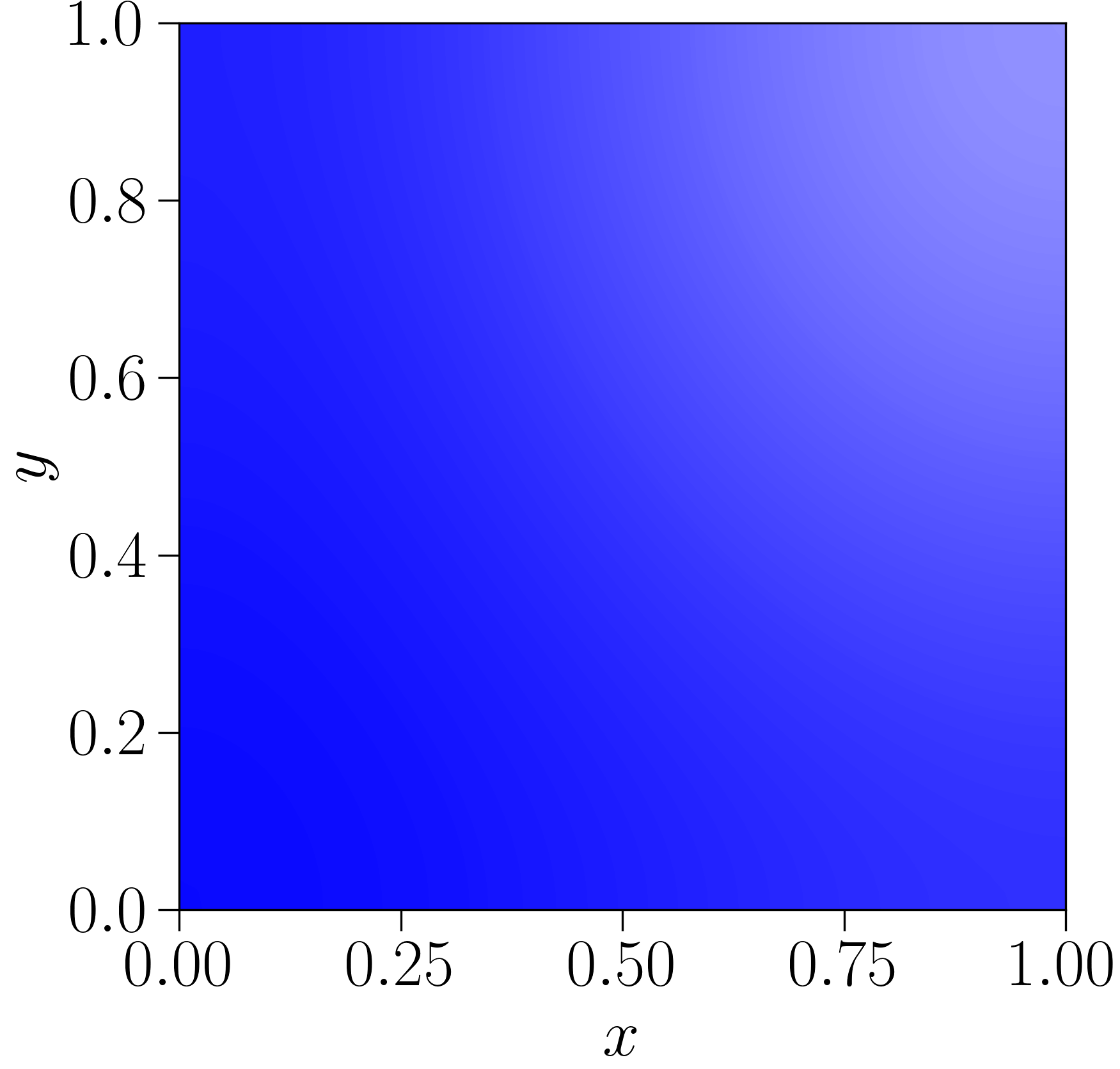}
			\end{minipage}
		}\\
		\subfloat{
			\begin{minipage}[b]{0.235\textwidth}
				\centering
				\includegraphics[width=1\linewidth]{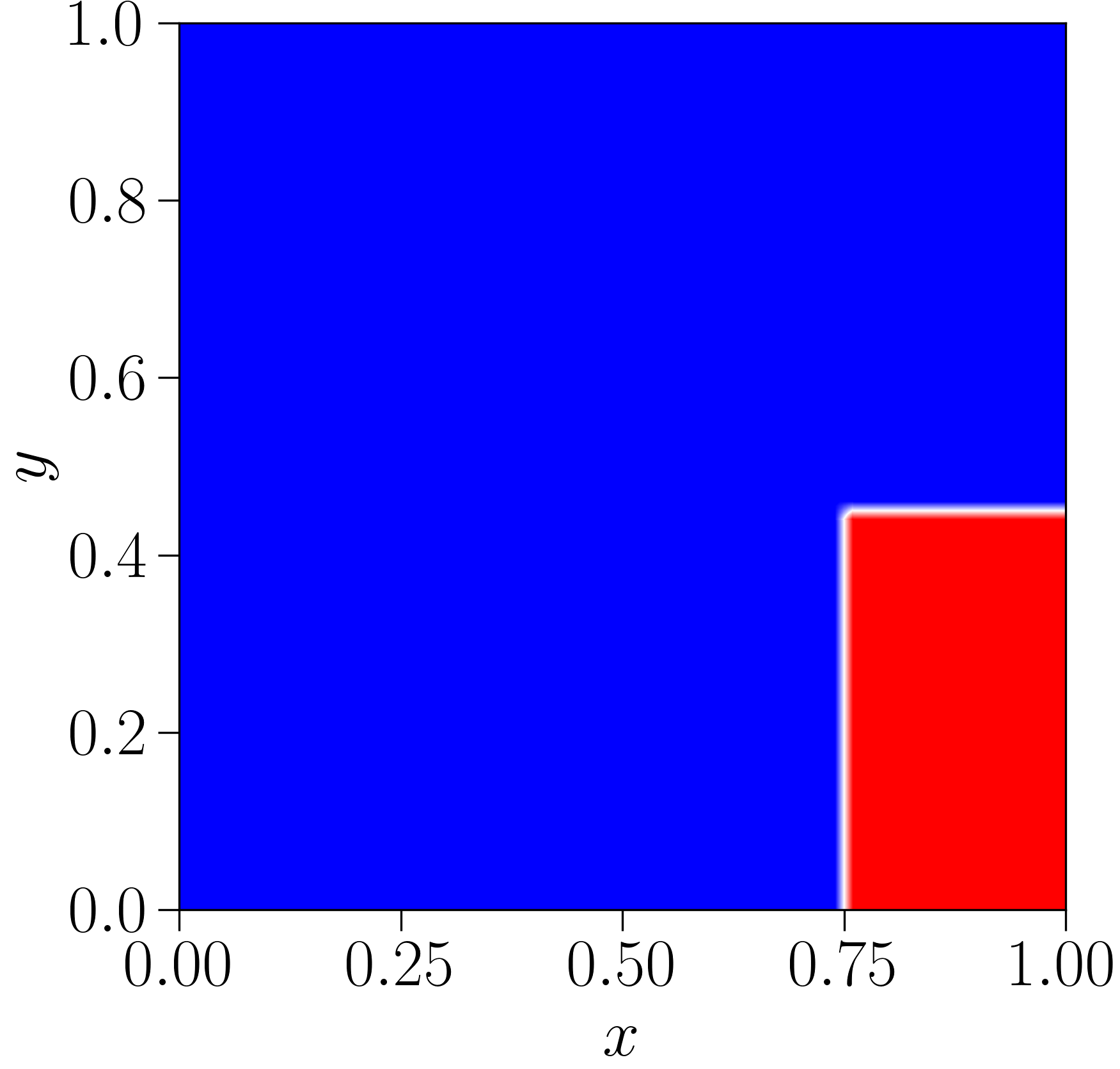}
			\end{minipage}
		}
		\subfloat{
			\begin{minipage}[b]{0.235\textwidth}
				\centering
				\includegraphics[width=1\linewidth]{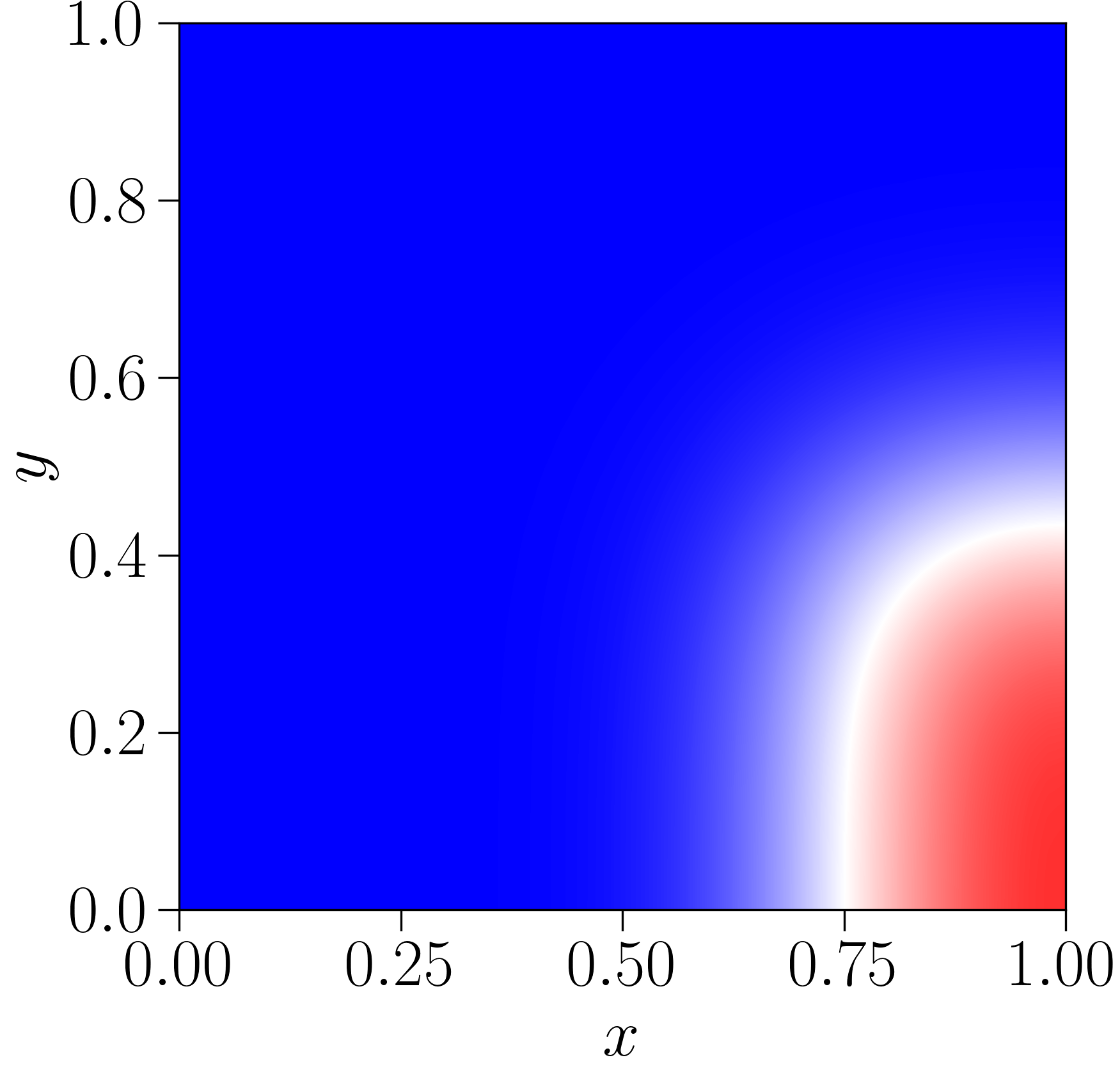}
			\end{minipage}
		}
		\subfloat{
			\begin{minipage}[b]{0.235\textwidth}
				\centering
				\includegraphics[width=1\linewidth]{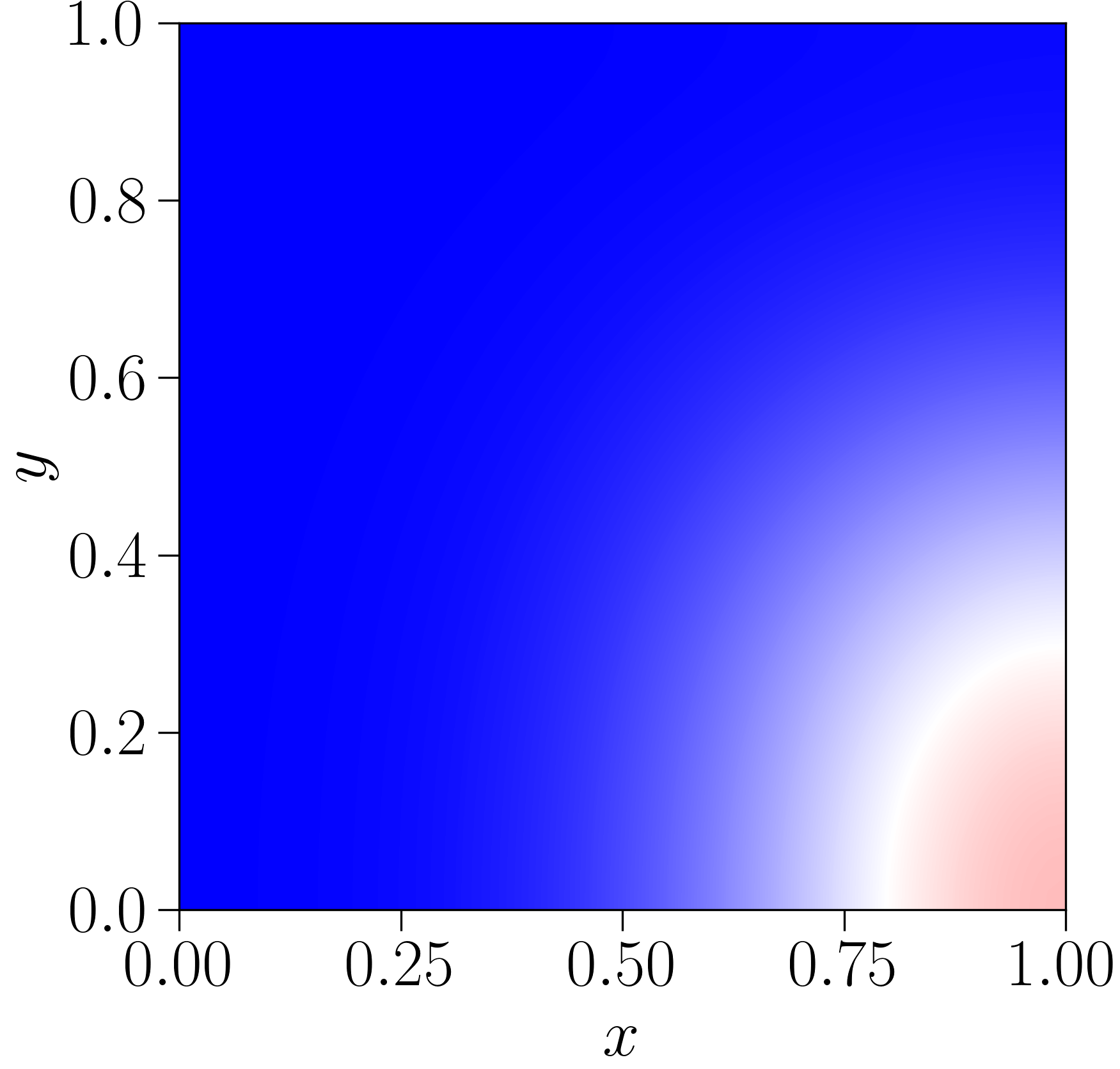}
			\end{minipage}
		}
		\subfloat{
			\begin{minipage}[b]{0.235\textwidth}
				\centering
				\includegraphics[width=1\linewidth]{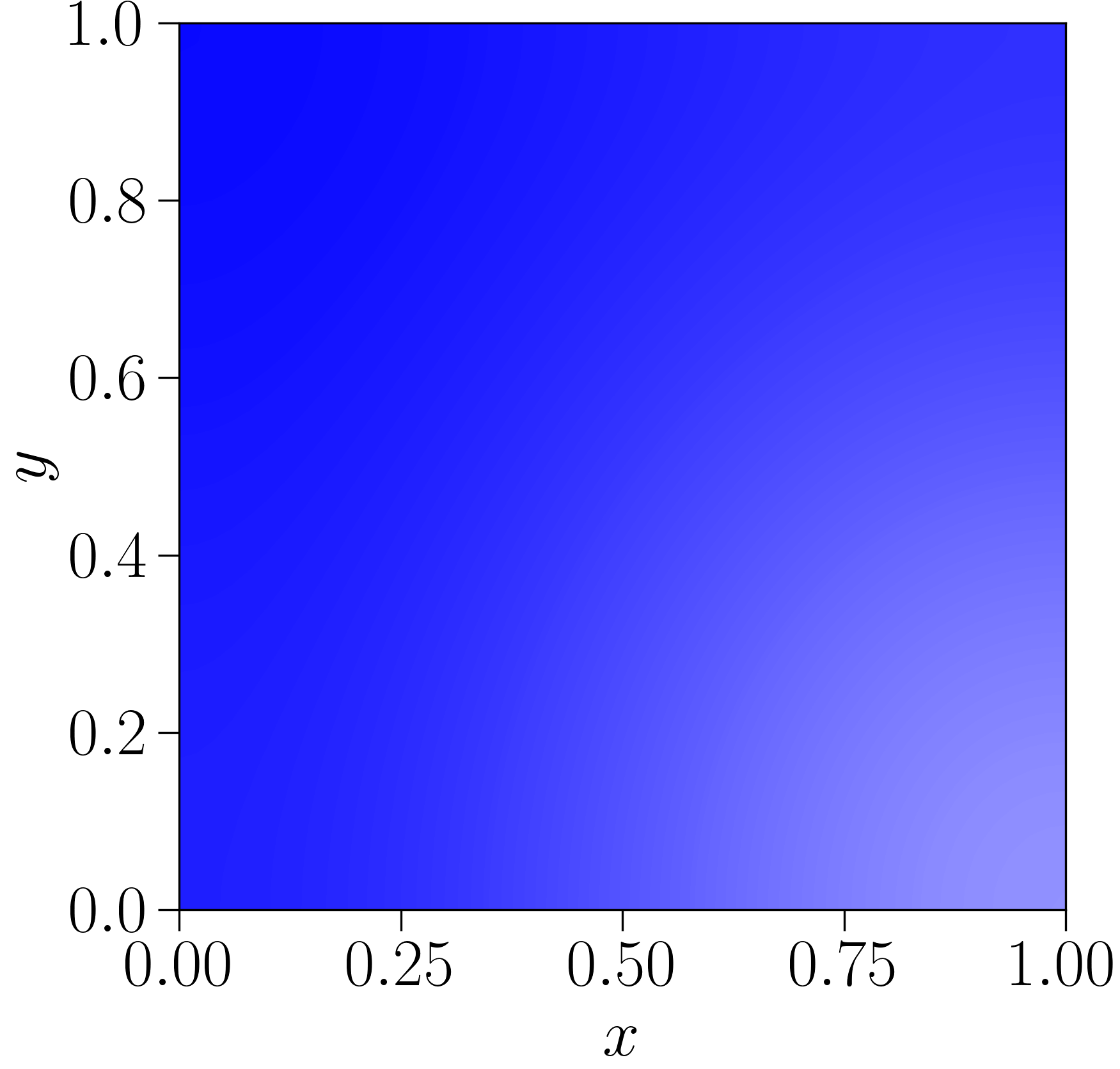}
			\end{minipage}
		}
	}
\end{minipage}
\begin{minipage}[c]{0.05\textwidth}%
	\centering
	\includegraphics[width=0.7\linewidth,]{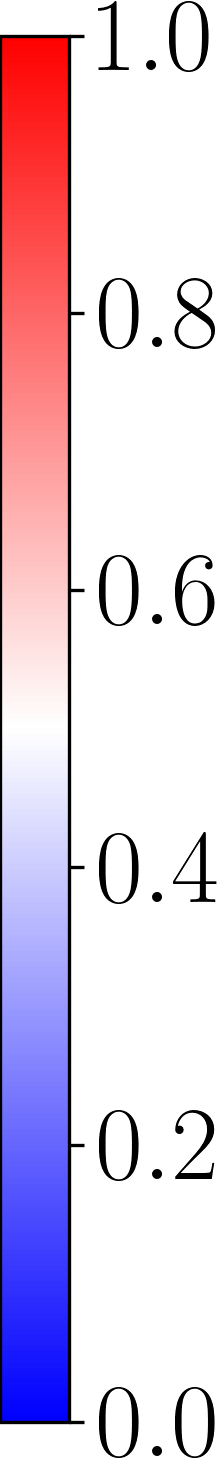}
\end{minipage}
\caption{Snapshots of $p^k_h$ (top) and $n^k_h$ (bottom) of the benchmark test of initial discontinuous concentrations at $t^k=0, 0.01, 0.03, 0.1$ (from left to right).}
\label{fig:IDC contour}
\end{figure}

\begin{figure}[h!]
\centering
\begin{minipage}[c]{1\textwidth}
	\centering{
		\subfloat{ 
			\begin{minipage}{0.35\textwidth}
				\centering
				\includegraphics[width=1\linewidth]{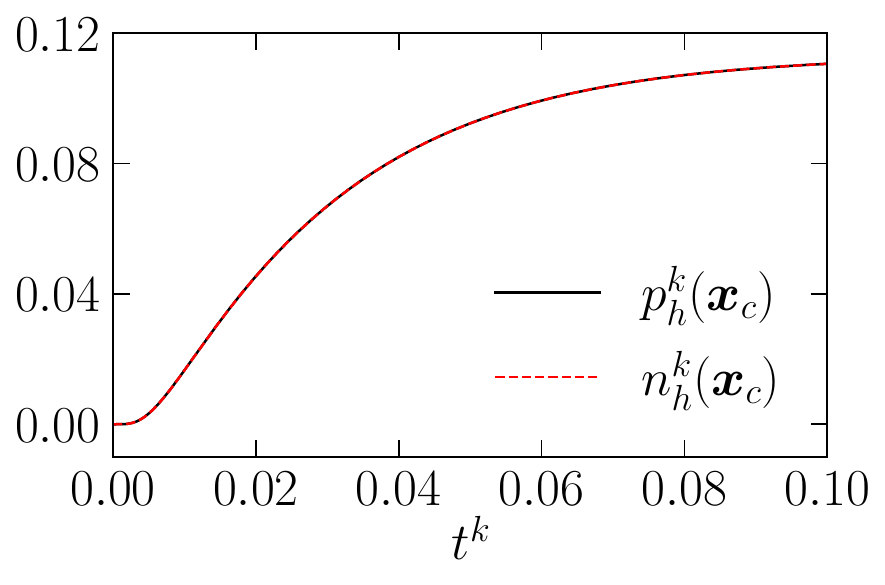}
			\end{minipage}
		}
		\subfloat{
			\begin{minipage}{0.7\textwidth}
				{
					\centering
					\begin{tabular}{ccccccc}
						\hline\bigstrut
						$t^k$   & 0.01 & 0.03 & 0.05 & 0.07 & 0.1  \\ \hline
						\bigstrut
						$p^k$ & 0.016555 & 0.066934 & 0.092285 & 0.103986 & 0.110558  \\
						$n^k$ & 0.016555 & 0.066934 & 0.092285 & 0.103986 & 0.110558  \\ \hline
					\end{tabular}
				}
			\end{minipage}
		}
	}
\end{minipage}%
\caption{Some results of the benchmark test of initial discontinuous concentrations. Left: $p_h^k(\boldsymbol{x}_c)$ and $n_h^k(\boldsymbol{x}_c)$ versus $t^k$. Right: $p_h^k(\boldsymbol{x}_c)$ and $n_h^k(\boldsymbol{x}_c)$ at some particular discrete time instants. $\boldsymbol{x}_c$ is the center of the computation domain, i.e. $\boldsymbol{x}_c = (0.5, 0.5) $.} 
\label{fig: IDC values}
\end{figure}

The last test is for a classic two-dimensional benchmark problem \cite{dehghan2023optimal} which has strong discontinuities in $p^0$ and $n^0$ within the computational domain $\Omega=(0, 1)^2$. The initial conditions are
\[
p^0 = \left\lbrace
\begin{aligned}
&1 \quad && \text{if } \boldsymbol{x}\in (0.75, 1) \times \left( \frac{11}{20}, 1\right) \\
&10^{-6} && \text{else}
\end{aligned}
\right.,
\]
\[
n^0 = \left\lbrace
\begin{aligned}
&1 \quad && \text{if } \boldsymbol{x}\in (0.75, 1) \times \left(0,  \frac{9}{20}\right) \\
&10^{-6} && \text{else}
\end{aligned}
\right.,
\]
$\boldsymbol{u}^0 = \boldsymbol{0}$ and $\omega^0 = 0$. And the boundary conditions are \eqref{Eq: boundary conditions}. The mesh used here is a uniform mesh of $50\times50$ square elements. The degree of finite-dimensional spaces is $N=1$. With this configuration, we can numerically compute $\psi_h^0$ by solving the discrete Poisson equation \eqref{eq sche c} at $t^0=0$. The simulation is then performed until $t^k=1$ using the present method at $\varDelta t = 10^{-4}$, and results are shown in Fig.~\ref{fig:IDC contour} and Fig.~\ref{fig: IDC values}. 

In Fig.~\ref{fig:IDC contour}, snapshots of $p^k_h$ and $n^k_h$ at $t^k\in\left\lbrace0,0.01,0.03,0.1\right\rbrace$ are present. It is seen that the initially locally distributed ions gradually pervade the whole domain. These results agree with those in the literature well\cite{dehghan2023optimal}. More quantitative results are given in Fig.~\ref{fig: IDC values}, where the values of $p^k_h $ and $n^k_h$ are in the center of the domain $\boldsymbol{x}_c=(0.5, 0.5)$. We can find that $p^k_h $ and $n^k_h$ have the same values (to a machine precision) at $\boldsymbol{x}_c$, which is consistent with the symmetric configuration of this problem. 

\section*{Acknowledgements}
 Yang is partially supported by the Science and Technology Project of Guangxi (Guike AD25069086) and is supported by the Natural Science Foundation of Guangxi (grant no. 2025GXNSFAA069683) and National Natural Science Foundation of China (grant no. 12561072). Zhang is partially supported by National Natural Science Foundation of Guangxi under grant no. 2024JJB110005. Peng is supported by the Innovation Project of GUET Graduate Education.


%

{\small
\bibliographystyle{elsarticle-num}
\bibliography{mybib}
}

\appendix
\section{} \label{App: A}
In this appendix, we provide the detailed derivation for the inequalities involving the logarithmic terms used in $E_1$ and
estimate the magnitude of the term$\left\|\nabla u_h^k - \nabla u_h^{k-1}\right\|^2$.
\begin{lemma}
	Let $f(t) = t \ln t - t$ be a convex function and define the functional $S[x] = \int_{\Omega} f(x) \mathrm{d}w$. For $x, y$ bounded such that $0 < m \le x, y \le M < \infty$, there exists a constant $C > 0$ such that:
	
	$$\langle x-y, \ln x \rangle - C\|x-y\|^2 \ge \langle f(x) - f(y), 1 \rangle$$
\end{lemma}
\paragraph{Proof}
Consider the second-order Taylor expansion of $f(t)$ at $x$. 
For any $x, y \in [m, M]$, there exists a value $\xi$ between $x$ and $y$ such that
\[f(y) = f(x) + f'(x)(y-x) + \frac{1}{2}f''(\xi)(y-x)^2.\]
Substituting $f'(t) = \ln t$ and $f''(t) = \frac{1}{t}$, we obtain:
\[f(y) = f(x) + (y-x) \ln x + \frac{1}{2\xi}(y-x)^2.\]
Rearranging the terms yields
\[(x-y) \ln x - \frac{1}{2\xi}(x-y)^2 = f(x) - f(y).\]
Given that $\xi \le M$, it follows that $\frac{1}{2\xi} \ge \frac{1}{2M}$. By setting $C' = \frac{1}{2M} > 0$, we have
\[(x-y) \ln x - C'(x-y)^2 \ge f(x) - f(y).\]
Integrating both sides over the domain $\Omega$ with respect to $w$ gives
\[\int_{\Omega} (x-y) \ln x \mathrm{d}w - C' \int_{\Omega} (x-y)^2 \mathrm{d}w \ge \int_{\Omega} (f(x) - f(y)) \mathrm{d}w.\]
In terms of the $L^2$ inner product and its norm, it is expressed as
\[\langle x-y, \ln x \rangle - C'\|x-y\|^2 \ge \langle f(x) - f(y), 1 \rangle.\]
Replacing $x$ and $y$ with $p_h^k$ and $p_h^{k-1}$ (or $n_h^k$ and $n_h^{k-1}$) completes the proof.

\begin{lemma} \label{lem:stability}
	Assume that  $p_h^k \geq m > 0$ for all $k$. Let $V_h$ be a finite element space satisfying the inverse inequality. Then, there exists a constant $C$ independent of the mesh size $h$ such that the following estimate holds:
	\begin{equation}
		\left\|\nabla \mu_h^k - \nabla \mu_h^{k-1}\right\|^2 \leq C h^{-2} \left( \left\| p_h^k - p_h^{k-1} \right\|^2 + \left\| \nabla \psi_h^k - \nabla \psi_h^{k-1} \right\|^2 \right).
	\end{equation}
\end{lemma}
\paragraph{Proof}
It follows from equation \eqref{eq scheme 2a} that:	
\[
\left\langle \mu_h^k , \xi_h \right\rangle = \left\langle \ln p_h^k  + \psi_h^k, \xi_h \right\rangle, \quad \forall \xi_h \in V_h.
\]
\[
\left\langle \mu_h^{k-1} , \xi_h \right\rangle = \left\langle \ln p_h^{k-1}  + \psi_h^{k-1}, \xi_h \right\rangle, \quad \forall \xi_h \in V_h.
\]
By considering the difference between two consecutive iterations, we obtain the identity:
\[
\langle \mu_h^k - \mu_h^{k-1}, \xi_h \rangle = \langle \ln p_h^k - \ln p_h^{k-1} + \psi_h^k - \psi_h^{k-1}, \xi_h \rangle, \quad \forall \xi_h \in V_h.
\]
Testing the above equation with $\xi_h = \mu_h^k - \mu_h^{k-1}$ and applying the Cauchy--Schwarz inequality, it follows that
\begin{equation}\label{eq:logest}
	\left\|\mu_h^k - \mu_h^{k-1}\right\| \leq \left\| \ln p_h^k - \ln p_h^{k-1} \right\| + \left\| \psi_h^k - \psi_h^{k-1} \right\|.
\end{equation}
To bound $\left\| \ln p_h^k - \ln p_h^{k-1} \right\|$, we note that for any $x \in \Omega$, the mean value theorem implies there exists a value $\zeta(x)$ between $p_h^k(x)$ and $p_h^{k-1}(x)$ such that
\begin{equation*}
	\left|\ln p_h^k(x) - \ln p_h^{k-1}(x)\right| = \left| \frac{1}{\zeta(x)} (p_h^k(x) - p_h^{k-1}(x)) \right|.
\end{equation*}
Given the uniform lower bound $p_h^k \geq m > 0$, it follows that $1/\zeta(x) \leq 1/m$. Squaring and integrating over $\Omega$ yields the $L^2$ estimate:
\[
\left\| \ln p_h^k - \ln p_h^{k-1} \right\|^2 \leq \frac{1}{m^2} \left\| p_h^k - p_h^{k-1} \right\|^2.
\]
Namely,
\begin{equation} \label{eq:log_est}
	\left\| \ln p_h^k - \ln p_h^{k-1} \right\| \leq \frac{1}{m} \left\| p_h^k - p_h^{k-1} \right\|.
\end{equation}
Substituting \eqref{eq:log_est} into \eqref{eq:logest}, we obtain the stability of the iteration in the $L^2$-norm:
\begin{equation} \label{eq:mu_L2_final}
	\left\|\mu_h^k - \mu_h^{k-1}\right\| \leq \frac{1}{m} \left\| p_h^k - p_h^{k-1} \right\| + \left\| \psi_h^k - \psi_h^{k-1} \right\|.
\end{equation}
Since $\mu_h^k - \mu_h^{k-1}$ belongs to the finite element space $V_h$, we can invoke the standard inverse inequality $\left\|\nabla v_h\right\| \leq C_{\mathrm{inv}} h^{-1} \left\|v_h\right\|$. Applying this to \eqref{eq:mu_L2_final} gives
\[
\left\|\nabla \mu_h^k - \nabla \mu_h^{k-1}\right\| \leq C_{\mathrm{inv}} h^{-1} \left( \frac{1}{m} \left\| p_h^k - p_h^{k-1} \right\| + \left\| \psi_h^k - \psi_h^{k-1} \right\| \right).
\]
Finally, by  inequality $(a+b)^2 \leq 2a^2 + 2b^2$ and the Poincaré inequality $\left\|\psi_h^k - \psi_h^{k-1}\right\| \leq C_P \left\|\nabla (\psi_h^k - \psi_h^{k-1})\right\|$, we conclude that
\[
\left\|\nabla \mu_h^k - \nabla \mu_h^{k-1}\right\|^2 \leq C h^{-2} \left\| p_h^k - p_h^{k-1} \right\|^2 + C h^{-2} \left\| \nabla \psi_h^k - \nabla \psi_h^{k-1} \right\|^2,
\]
where the constant $C$ depends on $m$, $C_{\mathrm{inv}}$, and $C_P$, but remains independent of the mesh size $h$.
The proof is complete.

\end{document}